\documentclass[11pt,a4paper]{article}
\usepackage{xypic}
\usepackage{french}
\usepackage{amssymb}
\usepackage{amsfonts}
\usepackage{amsmath}
\usepackage[mathscr]{eucal}
\usepackage{latexsym}
\begin{document}
\author{Gentiana Danila\\
Institut de Math\'ematiques de Jussieu, 
UMR 7586 du CNRS, Case Postale 7012\\
2, Place Jussieu, 
75251 Paris Cedex 05\\
e-mail: gentiana@math.jussieu.fr}
\title{R\'esultats sur la conjecture de dualit\'e \'etrange sur le
  plan projectif }
\date{1er avril 2000}
\maketitle
\newtheorem{theor}{Th\'eor\`eme}[section]
\newtheorem{prop}[theor]{Proposition}
\newtheorem{cor}[theor]{Corollaire}
\newtheorem{lemme}[theor]{Lemme}
\newtheorem{slemme}[theor]{Sous-lemme}
\newtheorem{defi}[theor]{D\'efinition}
\newtheorem{defiprop}[theor]{D\'efinition-Proposition}
\newtheorem{conj}[theor]{Conjecture}
\newtheorem{rem}[theor]{Remarque}
\newtheorem{rems}[theor]{Remarques}
\newtheorem{nota}[theor]{Notations}
\newtheorem{rappel}[theor]{Rappel}

\def\vs{\vskip 0.3cm}
\def\coker{{\rm coker}\, }
\def\im{{\rm im}\, }
\def\rmK{{\rm K}}
\def\H{{\rm H}}
\def\sl{{\rm SL}}
\def\gl{{\rm GL}}
\def\em{{\rm M}}
\def\es{{\rm S}}
\def\Hilb{{\rm Hilb}}
\def\Hom{{\rm Hom}}
\def\Pic{{\rm Pic}\,}
\def\Tor{{\rm Tor}\,}
\def\Ext{{\rm Ext}\,}
\def\uTor{\underline{{\rm Tor}}\,}
\def\uExt{\underline{{\rm Ext}}\,}
\def\uHom{\underline{{\rm Hom}}\,}
\def\uC{\underline{C}}
\def\underp{{\underline{p}}}
\def\tang{{\bf T}}
\def\div{{\rm div}\,}
\def\divs{{\rm div}\,\theta}
\def\deter{{\rm det}\, }
\def\deg{{\rm deg}\, }
\def\haut{{\rm ht}\,}
\def\profo{{\rm prof}\,}
\def\codim{{\rm codim}\, }
\def\sym{{\rm Sym\,}}
\def\supp{{\rm supp}\, } 
\def\Ker{{\rm Ker}\, }
\def\gr{{\rm gr}}
\def\cf{{\it cf.}}

\def\D{{\mathcal{D}}} 
\def\O{{\mathcal{O}}}
\def\R{{\mathcal{R}}}
\def\gotu{{\mathfrak{u}}}
\def\gothas{{\mathfrak{h}}}
\def\V{{\mathcal V}}
\def\maG{{\mathcal{G}}}
\def\maF{{\mathcal{F}}}
\def\maA{{\mathcal{A}}}
\def\maB{{\mathcal{B}}}
\def\maC{{\mathcal{C}}}
\def\maE{{\mathcal{E}}}
\def\maL{{\mathcal{L}}}
\def\maP{{\mathcal{P}}}
\def\maT{{\mathcal{T}}}
\def\maH{{\mathcal{H}}}
\def\maI{{\mathcal{I}}}
\def\sigm{{\mathfrak{S}}_m}
\def\maQ{{\mathcal{Q}}}
\def\maR{{\mathcal{R}}}
\def\has{{\mathcal H}}
\def\K{{\mathcal K}}

\def\xim{\Xi_{m}}
\def\kp{{\rm K}(\pp)}
\def\proj{{\mathbb{P}}}
\def\pp{\proj_2}
\def\pps{\pp^*}
\def\hil{\Hilb^m(\pp)}
\def\hilx{X^{\mbox{}^{[m]}}}
\def\hilxp{X^{\mbox{}^{[m+1]}}}
\def\hilxx{\hilx\times X}
\def\hilxmm{X^{\mbox{}^{[m,m+1]}}}
\def\en{{\mathbb N}}
\def\zed{{\mathbb Z}}
\def\comp{{\mathbb C}}
\def\qiu{{\mathbb Q}}
\def\xm{X^m}
\def\smx{\es^m(X)}
\def\emce{{\em_c}}
\def\emen{{\em_n}}
\def\emu{{\em_u}}
\def\emceu{\emce\times\emu}
\def\emd{{\em_{d{\mathfrak{u}}}}}
\def\emunu{{\em_{\mathfrak{u}}}}
\def\emdoi{{\em_{2{\mathfrak{u}}}}}
\def\emtrei{{\em_{3{\mathfrak{u}}}}}
\def\zu{Z_{univ}}
\def\zup{Z^{\prime}_{univ}}
\def\zupp{Z^{\prime\prime}_{univ}}
\def\xiu{\xi_{univ}}
\def\xipu{\xi^{\prime}_{univ}}
\def\xippu{\xi^{\prime\prime}_{univ}}
\def\etau{\eta_{univ}}
\def\etapu{\eta^{\prime}_{univ}}
\def\etappu{\eta^{\prime\prime}_{univ}}
\def\blup{{\rm{Bl}}_{\Xi_m}(\hilx\times X)}
\def\ix{{{I}}_{\xim}}
\def\oxm{\O_{\xim}}
\def\pix{\proj(\ix)}
\def\emceu{\emce\times\emu}
\def\Omegace{\Omega_c}
\def\Omegau{\Omega_u}
\def\Omegaceu{\Omega_c\times\Omegau}
\def\deomega{\D_{\Omega}}
\def\sigmaom{\sigma_{\Omega}}
\def\gece{G_c}
\def\geu{G_u}
\def\geceu{\gece\times\geu}
\def\efce{\maF_c}
\def\efu{\maF_u}
\def\ghec{\maG_c}
\def\gheu{\maG_u}
\def\esceesu{(s_c,s_u)}
\def\geesce{G_{s_c}}
\def\geesu{G_{s_u}}
\def\geesceu{\geesce\times\geesu}
\def\deesceu{\D_{\esceesu}}
\def\efesce{\maF_{s_c}}
\def\gheesu{\maG_{s_u}}
\def\deceesu{\D_{c,s_u}}
\def\emceesu{\emce\times\{s_u\}}
\def\esceemu{\{s_u\}\times\emu}
\def\emcepedoi{\em_{\frac{c}{2}}}

\def\vmu{\O(-1)^{\mbox{}^{[m]}}}
\def\decu{\D_{c,u}}
\def\dees{\D_S}
\def\sigmaes{\sigma_S}
\def\de{{\mathfrak{d}}}
\def\dece{\D_c}
\def\deu{\D_u}
\def\deuce{\deu\boxtimes \dece}
\def\deunu{{\D_{\mathfrak{u}}}}
\def\dedoi{{\em_{2{\mathfrak{u}}}}}
\def\detrei{{\em_{3{\mathfrak{u}}}}}
\def\dede{{\D^{\otimes d}}}
\def\deka{{\D^{\otimes k}}}
\def\deji{{\D^{\otimes j}}}
\def\demi{{\D^{-1}}}
\def\dei{{\D^{\otimes -i}}}
\def\dekai{{\D^{\otimes k-i}}}
\def\dekam{{\D^{\otimes k-m}}}
\def\dekaunu{{\D^{\otimes k-1}}}
\def\dejide{{\D^{\otimes -j-3\delta}}}
\def\dekaide{{\D^{\otimes -k+i-3\delta}}}
\def\dekaji{{\D^{\otimes m-k-3\delta-j}}}
\def\ded{{\D_{d{\mathfrak{u}}}}}
\def\vk{\O(k)^{\mbox{}^{[m]}}}
\def\vtrei{\O(3)^{\mbox{}^{[m]}}}
\def\vunu{\O(1)^{\mbox{}^{[m]}}}
\def\vkde{\vk\otimes\de}
\def\lm{L^{\mbox{}^{[m]}}}
\def\deA{{\mathfrak{d}}^A_m}
\def\ox{\omega_X}
\def\omegau{\omega_{\emu}}
\def\omegade{\omega_{\emd}}
\def\omegace{\omega_{\emce}}

\def\surto{\twoheadrightarrow}
\def\ra{\rightarrow}

\def\sigmaceu{\sigma_{c,u}}
\def\deceu{{\rm D}_{c,u}}
\def\decedeu{{\rm D}_{c,d{\mathfrak{u}}}}
\def\deceunu{{\rm D}_{c,{\mathfrak{u}}}}

\def\tens{\otimes}
\def\tna{\widetilde{\nabla}}
\def\tD{\widetilde{D}}
\def\tnu{\widetilde{\nu}}


\def\fs{{faisceau }}
\def\fx{{faisceaux }}
\def\alg{{alg\'ebrique }}
\def\algs{{alg\'ebriques }}
\def\th{{th\'eor\`eme }}
\def\iso{{isomorphisme }}
\def\rep{{repr\'esentation }}
\def\reps{{repr\'esentations }}
\def\irr{{irr\'eductible }}
\def\irrs{{irr\'eductibles }}
\def\fib{{fibr\'e }}
\def\fibs{{fibr\'es }}
\def\mor{{morphisme }}
\def\mors{{morphismes }}
\def\sur{{surjectif }}
\def\co{{coh\'erent }}
\def\cohs{{coh\'erents }}
\def\app{{application }}
\def\apps{{applications }}
\def\cep{{caract\'eristique d'Euler-Poincar\'e }}
\def\Gro{{Grothendieck }}
\def\espmod{{espace de modules }}
\def\espmods{{espaces de modules }}
\def\ss{{semi-stable }}
\def\sss{{semi-stables }}
\def\resp{{respectivement }}
\def\inve{{inversible }}
\def\inves{{inversibles }}
\def\fibd{{\fib d\'eterminant }}
\def\sct{{section }}
\def\scts{{sections }}
\def\cano{{canonique }}
\def\canos{{canoniques }}
\def\lin{{lin\'eaire }}
\def\lins{{lin\'eaires }}
\def\schil{{sch\'ema de Hilbert }}

{\bf Abstract~:} {\small Le Potier's ``Strange Duality'' conjecture
  gives an isomorphism between the space of sections of the
  determinant bundle on two different moduli spaces of semi-stable sheaves on the
  projective plane $\pp$. If we consider two orthogonal classes $c,u$ in the
  Grothendieck algebra $\kp$ such that $c$ is of positive rank and $u$
  of rank zero,  we call $\emce$ and $\emu$ the moduli spaces of
  semi-stable sheaves of class $c$, respectively $u$ on $\pp$. There
  exists on $\emce$ (resp. $\emu$) a
  determinant bundle $\deu$ (resp. $\dece$) and the product fibre
  bundle $\deuce$ on the product space $\emceu$ has a canonical
  section $\sigmaceu$ which provides a linear application
  $\deceu:\H^0(\emu,\dece)^*\to\H^0(\emce,\deu)$. If $\emce$ is not
  empty, $\deceu$ is conjectured to be an isomorphism. We prove the
  conjecture in the particular case where $c$ is of rank $2$, 
  zero first Chern class and second Chern class $c_2(c)\le 5$, and $u$ is of
  degree $d(u)\le 3$ and zero Euler-Poincar\'e characteristic. In
  addition we give the generating series $P(t)=\sum_{k\ge
  0}t^kh^0(\emce,\deu^{\otimes k})$ for $c_2(c)=3$, $c_2(c)=4$,
  $d(u)=1$, for
  the particular classes $c$ and $u$ considered above.}



Running heads: R\'esultats sur la dualit\'e \'etrange sur le
  plan projectif

\section{Introduction}

\label{sec:int}
La motivation principale de cet article est de fournir des exemples en
faveur de la dualit\'e \'etrange sur le plan projectif conjectur\'ee
par Le Potier (\cite{LeP-exposes}). On
consid\`ere l'alg\`ebre de Grothendieck $\kp$ des classes de \fx \algs
\cohs sur $\pp$. C'est un groupe ab\'elien  isomorphe \`a
$\zed^3$, un \iso \'etant donn\'e par le rang, la classe de Chern et la \cep (ceci nous permet de d\'esigner chaque classe $c\in\kp$
par le triplet form\'e par son rang $r$, sa premi\`ere classe de Chern
$c_1$ et sa \cep $\chi$, ou bien, lorsque c'est indiqu\'e, sa deuxi\`eme
classe de Chern $c_2$). Elle est munie d'une multiplication et d'une
forme bilin\'eaire donn\'ee par $ \langle c,u \rangle=\chi(c\cdot u)$. Pour
deux classes $c$ et $u$ orthogonales, de rang $r>0$ et respectivement
$0$, on note $\emce$ et $\emu$ les \espmods des \fx \sss sur $\pp$ de
classe $c$ et respectivement $u$. Sur chacun de ces espaces il existe
un \fib \inve $\deu$ et \resp $ \dece$, appel\'e \fib d\'eterminant. Alors le \fib
produit tensoriel externe $\deuce$ sur  $\emceu$ a une \sct \cano
$\sigmaceu$, qui fournit une \app \lin
$$\deceu:\H^0(\emu,\dece)^*\to\H^0(\emce,\deu)$$
appel\'ee \mor de dualit\'e \'etrange. Remarquons que le groupe $\sl(3)$
agit sur $\pp$. Il agit ainsi sur les espaces de modules $\emce,
\emu$, et sur les \fibs d\'eterminants $\dece, \deu$. Le \mor $\deceu$
est un \mor de $\sl(3)$-repr\'esentations.

\vs{\bf Conjecture } (J. Le Potier) {\it Si $\em_c$ est non-vide alors  
le morphisme de dualit\'e \'etrange est un isomorphisme.}
\vs  On  consacre la section \ref{section1} \`a la
construction et l'interpr\'etation g\'eom\'etrique du \mor $\deceu$. On va se
restreindre ici au cas $c=(2,0,2-n)$ (qui recouvre le cas o\`u la
premi\`ere classe de Chern est paire), et $u=d(0,1,0)$. Le cas $d=1$
et $n\le 19$ a \'et\'e analys\'e dans l'article \cite{D1}. Le
r\'esultat principal est:
\begin{theor}
\label{dual}
Soit $\gotu=(1,0,0)\in\kp$. Si $r=2, c_1=0$ et $n=c_2\le 5$ (i.e. $c=(2,0,2-n)$) et $u=(0,d,0)=d\gotu$, alors l'application lin\'eaire
$$D_{c,d\gotu}:\H^0(\em_{d\gotu},\dece)^*\to\H^0(\em_c,\D_{\gotu}^{\otimes
  d})$$
est un isomorphisme pour $d=2,3$, c'est-\`a-dire que dans ces
  conditions la conjecture de dualit\'e \'etrange est vraie. 
\end{theor}

Au paragraphe \ref{lesespaces} on utilise \cite{LeP2} pour d\'ecrire les
espaces de modules $\emd$. Il existe un \mor $\pi:\emd\to C_d$
(espaces des courbes de degr\'e $d$ dans $\pp$) qui associe au \fs $G$
l'\'equation de son support sch\'ematique. C'est un \iso pour $d=1,2$ et
un \mor dont la fibre g\'en\'erique est de dimension $1$ pour $d=3$. Ceci
nous permet de calculer $\H^0(\emd,\dece)$. 

Au paragraphe \ref{inject} on d\'emontre l'injectivit\'e du \mor $\deceu$ en
utilisant l'interpr\'etation g\'eom\'etrique du \th
\ref{pro:cadru} (iv). On utilise les propri\'et\'es de $\emd$
\'etablies dans la section \ref{lesespaces}.

On calcule au paragraphe \ref{1}
les espaces $\H^0(\em_c,\D_{\gotu}^{\otimes
  2})$ et $\H^0(\em_c,\D_{\gotu}^{\otimes
  3})$ en tant que $\sl(3)$-re\-pr\'e\-sen\-ta\-tions, selon la technique
d\'evelopp\'ee dans l'article \cite{D1}.

\begin{prop}
\label{4.0.7}
Avec les notations du \th pr\'ec\'edent, le $\sl(3)$-module $\H^0(\em_c,\D_{\gotu}^{\otimes
  2})$ est isomorphe \`a $\es^n(\es^2E)$ et le $\sl(3)$-module $\H^0(\em_c,\D_{\gotu}^{\otimes
  3})$ est isomorphe \`a $\es^n(\es^3E)\oplus\es^{n-2}(\es^3E)$(o\`u $E=\H^0(\pp,\O(1))$ est la
  repr\'esentation standard de $\sl(3)$ ).
\end{prop}

Ceci nous permet de conclure la preuve du th\'eor\`eme.
De plus, nous calculons au paragraphe \ref{3} 
la dimension des espaces de sections de $\D_{\gotu}^{\otimes k}$ pour
$n=3,4$, qu'on \'ecrit sous forme de s\'erie de Poincar\'e~: 

\begin{theor}
\label{203204}

i) Pour l'espace de modules $\em_{(2,0,-1)}$ des \fx stables de rang
$2$ et classes de Chern $(0,3)$, la s\'erie de Poincar\'e de $\D_{\gotu}$,
$P(t)=\sum_{k\ge 0}t^kh^0(\D_{\gotu}^{\otimes k})$ est donn\'ee par
$$P(t)=\frac{1+t^2+t^4}{(1-t)^{10}}.$$

ii) Pour l'espace de modules $\em_{(2,0,-2)}$ des \fx semi-stables de rang
$2$ et classes de Chern $(0,4)$, la s\'erie de Poincar\'e de $\D_{\gotu}$  est donn\'ee par
$$P(t)=\frac{1+t+7t^2+7t^3+22t^4+7t^5+7t^6+t^7+t^8}{(1-t)^{14}}.$$

\end{theor}

\section{Morphisme de dualit\'e \'etrange}
\label{section1}
L'objet de cette partie est de pr\'esenter la conjecture de Le Potier
sur la dualit\'e \'etrange.

\subsection{L'alg\`ebre de Grothendieck $\kp$}

Si $S$  est une vari\'et\'e alg\'ebrique, on d\'esigne par $\rmK(S)$ 
le groupe de Grothendieck  des classes de faisceaux alg\'ebriques
coh\'erents sur $S$.
Pour un \fs $F$ on
note $[F]$ sa classe dans le groupe $\rmK(S)$. Dans ce qui suit, on aura \`a
consid\'erer en particulier  le groupe de Grothendieck $K(\pp)$:  c'est
un groupe ab\'elien libre de rang 3; 
l'application
$$\begin{array}{ccccc}
\phi&:&\kp&\to&\zed^3\\
&&[F]&\mapsto&(rg(F),c_1(F),\chi(F))
\end{array}$$
qui \`a la classe d'un \fs $F$ associe le rang $r$ de $F$, la classe de
Chern $c_1$ de $F$, et la
caract\'eristique d'Euler-Poincar\'e $\chi$ de $F$, est un \iso de groupes ab\'eliens. On notera un \'el\'ement de $\kp$ par son
image $(r,c_1,\chi)$ dans $\zed^3$. Si $S$ est lisse, une loi de  multiplication sur $\rmK(S)$ est d\'efinie  en prolongeant par
lin\'earit\'e la loi de multiplication d\'efinie pour $F$ et $G$ \fx \algs
coh\'erents sur $S$ par
$$F\cdot G=\sum_p(-1)^p\underline{Tor}_p(F,G)$$

Ce produit se r\'eduit au produit tensoriel usuel si  l'un des deux \fx $F$ ou $G$
est localement libre. 

On parle alors d'alg\`ebre de Grothendieck.

On note $\eta=[\O_l]$  la classe du \fs structural
d'une droite $l$, et $\eta^2=[\O_p]$ celle du \fs structural d'un
point. En tant qu'alg\`ebre, $\kp$ est isomorphe \`a
$\zed[\eta]/(\eta^3)$. On munit $\kp$ de la forme bilin\'eaire donn\'ee par $\langle c,u\rangle =\chi(c\cdot
u)$. 
Dans la suite l'orthogonalit\'e sera prise relativement \`a cette forme.
On a aussi sur
$K(\pp)$  une involution $u\mapsto u^*$  qui associe \`a la classe d'un
fibr\'e vectoriel celle de son dual.

\subsection{Espaces de modules et fibr\'es d\'eterminants}
 
Soit $S$ une vari\'et\'e. On note $pr_1$ la projection $S\times \pp\to
S$ et $pr_2$ la projection $S\times \pp\to\pp$. Pour un \fs $\maF$ sur $S\times \pp$, on note
$$pr_{1!}([\maF]):=[R^0pr_{1*}\maF]-[R^1pr_{1*}\maF]+[R^2pr_{1*}\maF]$$
dans le groupe $\rmK(S)$. Cela d\'efinit, par lin\'earit\'e, une application
$pr_{1!}:\rmK(S\times\pp)\to\rmK(S)$. 

Si $F$ est un \fs \co sur $S$ qui
admet une r\'esolution finie par des \fx localement libres $A_i$~:
$$0\to A_n\to A_{n-1}\to\cdots\to A_0\to F\to 0$$
on introduit le
\fs \inve
$$\det F=\det A_0\otimes(\det A_1)^{-1}\otimes\cdots\otimes(\det
A_n)^{(-1)^n}.$$
L'application qui \`a $F$ associe son d\'eterminant $\det F$
est multiplicative sur les suites exactes. 

Soient $c\in\kp$ une classe de Grothendieck de rang $r>0$ et $\em_c$
l'espace de modules des \fx semi-stables de classe de Grothendieck
$c$. C'est une vari\'et\'e \alg projective irr\'eductible normale, de
dimension $D=1-\langle~c^*,c\rangle$ o\`u $c^*$ est la classe duale
de $c$. On note $c^{\bot}$ le sous-espace de $\kp$ des classes
orthogonales \`a $c$. Dans
\cite{Dre}, Dr\'ezet a construit un morphisme surjectif de groupes
$\lambda_c:c^{\bot}\to\Pic(\em_c)$ caract\'eris\'e par la propri\'et\'e  universelle
suivante~:

Pour toute famille plate $\maF$ de \fx semi-stables de classe de
Grothendieck $c$, param\'etr\'ee par une vari\'et\'e
alg\'ebrique $S$, la classe $\lambda_{\maF}(u)=\det pr_{1!}(\maF\cdot
pr_2^*(u))$
d\'efinit un  fibr\'e inversible sur la vari\'et\'e $S$. Si
$f_{\maF}:S\to\em_c$ est le morphisme modulaire associ\'e \`a la famille
$\maF$, on a
$$f_{\maF}^*(\lambda_c(u))=\lambda_{\maF}(u)$$
et le fibr\'e  $\lambda_c(u)$ est le seul \`a \iso pr\`es qui satisfait \`a cette propri\'et\'e
pour toute famille plate $\maF$.

S'il existe un \fs universel $\maF$ sur ${\emce}$, d'apr\`es la
propri\'et\'e universelle  de $\lambda_{c}(u)$, il r\'esulte que
$\lambda_{c}(u)=\lambda_{\maF}(u)$. En g\'en\'eral, on prouve l'existence
de  $\lambda_{c}(u)$ en \'ecrivant $\emce=\Omega^{ss}/\sl(H)$ comme
quotient d'un ouvert dans un \schil par l'action d'un groupe r\'eductif, en
consid\'erant une famille universelle sur $\Omega^{ss}$ et en utilisant
un argument de descente.

\vs On note $\rmK_c$ le sous-$\zed$-module libre de rang $1$ de $c^{\bot}$
des classes de rang $0$. 
On appelle $\gotu$ le g\'en\'erateur positif de $\rmK_c$ (i.e. $c_1(\gotu)$ est un
multiple positif de la classe hyperplane $h$ dans $\H^2(\pp,\zed)$). On a
$\gotu=(0,\frac{r^c}{\delta},-\frac{c_1^c}{\delta})$  o\`u
$\delta={\rm pgcd}\,(r^c,c_1^c)$.

Le \fib $\D=\D_{\gotu}=\lambda_c(-\gotu)$ s'appelle \fibd de Donaldson
sur $\em_c$. Pour toute classe $u\in\rmK_c$ on introduit plus
g\'en\'eralement le fibr\'e inversible $\D_u=\lambda_c(-u)$; c'est donc un
multiple du fibr\'e d\'eterminant de Donaldson. 

Le probl\`eme du calcul de 
la dimension de l'espace de
sections $\H^0(\em_c,\D_u)$ a conduit 
Le Potier (\cite{LeP2}) \`a introduire $\em_u$, l'espace de modules des \fx
semi-stables de dimension $1$ de classe de Grothendieck $u$. La
notion de semi-stabilit\'e (resp. stabilit\'e) se g\'en\'eralise
(\cf \cite{LeP2}) pour
les \fx \algs
\cohs $F$ de dimension $1$.  On sait que $\em_u$ est encore
une vari\'et\'e \alg  irr\'eductible normale et que la classe $c$
(appartenant \`a $u^{\bot}$) permet de construire un fibr\'e inversible
$\D_c=\lambda_u(-c)$ sur $\em_u$, de la m{\^e}me mani\`ere que sur $\em_c$. Dans certains
cas l'espace $\em_u$ et le fibr\'e $\D_c$ sont plus faciles
\`a d\'ecrire.

\subsection{Construction du morphisme de dualit\'e \'etrange}

\label{morphisme}

On consid\`ere deux classes  $c,u\in\kp$ dans l'alg\`ebre de \Gro. On
suppose qu'elles sont orthogonales, que $r(c)>0$, que $r(u)=0$ et
$c_1(u)>0$. Le \mor de dualit\'e \'etrange est une cons\'equence de la
construction simultan\'ee des \fibs d\'eterminants $\deu$ sur $\emce$ et
$\dece$ sur $\emu$. Sa construction et son  interpr\'etation
g\'eom\'etrique sont r\'esum\'ees dans le \th suivant~:
\begin{theor}
  \label{pro:cadru}

i) Il existe une section canonique, d\'efinie \`a une constante pr\`es,
$\sigmaceu\in\H^0(\emceu,\deuce)$ qui s'annule exactement aux points $([F],[G])$ tels que $h^0(\pp,
F\tens G)=h^1(\pp,F\tens G)\ne 0$.

ii) La section $\sigmaceu$ d\'efinit une application lin\'eaire
$$\deceu:\H^0(\emu,\dece)^*\to\H^0(\emce,\deu).$$

iii) On note $\sigma_F$ la restriction de $\sigmaceu$ \`a
 $\{[F]\}\times\emu$. 
Si $\sigmaceu$ n'est pas identiquement nulle, l'association $F\mapsto[\sigma_F]$ d\'efinit  une application rationnelle
$$\Phi:\emce\to\proj\H^0(\emu,\dece).$$
 Si en outre pour tout $[F]\in\emce$, $\sigma_F$ n'est pas identiquement nulle, l'application $\Phi$ est r\'eguli\`ere.

iv) Si l'image du \mor $\Phi$
n'est pas contenue dans un hyperplan l'application $\deceu$ est injective.
\end{theor}

La conjecture de Le Potier est alors~:

\begin{conj}
\label{1.3.2}
Si $\em_c$ est non-vide alors  
le morphisme de dualit\'e \'etrange $\deceu$ est un isomorphisme.
\end{conj}

\par{\bf Preuve du \th \ref{pro:cadru}~:}

On commence par rappeler les r\'esultats suivants~:

\begin{lemme}
\label{lemme1}
  Soit $S$ une vari\'et\'e alg\'ebrique. Soient $\maF$ et $\maG$ des
  familles plates de \fx \sss sur $\pp$ de classes de \Gro $c$ et
  \resp $u$, param\'etr\'ees par $S$. Alors~:

a) Le \fs $\maG$ a une r\'esolution 
\begin{equation}
  \label{49}
 0\to \maQ\to\maR\to\maG\to 0 
\end{equation}
sur $S\times\pp$ par des \fx localement libres $\maQ$ et $\maR$.

b) Le \fs $\maF$ a une r\'esolution 
\begin{equation}
  \label{50}
 0\to \maA\to\maB\to\maF\to 0 
\end{equation}
sur $S\times\pp$ par des \fx localement libres $\maA$ et $\maB$. En
plus, on peut choisir $\maB$ tel que 
\begin{equation}
  \label{505}
  h^0(\maB_s\tens\maG_s)=0 {\rm \ pour \ tout \ } s\in S. 
\end{equation}

c) $\uTor_i(\maF,\maG)=0$ pour $i>0$.
\end{lemme}

\par{\bf Preuve~:}

La r\'esolution pour $\maG$ r\'esulte du fait que $\maG$ est pur de
dimension $1$ sur chaque fibre. La r\'esolution pour $\maF$ r\'esulte du
fait que la restriction de $\maF$ \`a chaque fibre est un \fs sans torsion sur
$\pp$. On peut changer le \fs $\maB$ dans la r\'esolution de $\maF$ en
$[\maB\tens pr_2^*\O(-n)]^m$, pour $n\ge 0$ et $m$ assez grand. Pour un choix de $n$ assez
grand on obtient 
$h^0(\maB_s\tens\maG_s)=0$ pour tout $s$. Il r\'esulte d'apr\`es a), b),
que $\uTor_i(\maF,\maG)=0$ pour $i\ge 2$ et que $\uTor_1(\maF,\maG)$
est inclus dans $\maG\tens\maA$ et dans $\maF\tens\maQ$. Comme  $\maG\tens\maA$
 est de torsion, et $\maF\tens\maQ$ sans torsion, on en d\'eduit que 
$\uTor_1(\maF,\maG)=0$. $\Box$

\vs Soient $S$, $\maF$, $\maG$ comme dans le lemme. D'apr\`es c) on a une
suite exacte courte
\begin{equation}
  \label{51}
  0\to \maA\tens\maG \stackrel{a}{\to}\maB\tens\maG \to
  \maF\tens\maG\to 0.
\end{equation}
On consid\`ere son image directe par $pr_{1*}$. Le lemme \ref{lemme1} b)
conduit \`a $pr_{1*}(\maA\tens\maG)=pr_{1*}(\maB\tens\maG)=0$. Le \fs
$\maG$ est de dimension $1$ dans les fibres, donc $R^2pr_{1*}(\maA\tens\maG)=R^2pr_{1*}(\maB\tens\maG)=0$.
Alors la suite
\begin{equation}
  \label{52}
  0\to pr_{1*}(\maF\tens\maG)\to
  R^1pr_{1*}(\maA\tens\maG)\stackrel{a}{\to}R^1pr_{1*}(\maB\tens\maG)\to
  R^1pr_{1*}(\maF\tens\maG)\to 0
\end{equation}
est exacte.

\begin{lemme}
  \label{lemme2}
Les \fx $R^1pr_{1*}(\maA\tens\maG)$ et $R^1pr_{1*}(\maB\tens\maG)$
sont localement libres de m\^eme rang sur $S$.
\end{lemme}

\par{\bf Preuve~:}

Puisque les familles $\maF$ et $\maG$ sont $S$-plates, on obtient que
$\maA\tens\maG$, $\maB\tens\maG$ et $\maF\tens\maG$ sont des familles
$S$-plates. Alors les suites (\ref{49}), (\ref{51}), (\ref{52}) sont
compatibles avec les changements de base $S'\to S$. En particulier,
pour $S'=\{s\}\in S$, on obtient  \`a partir de la suite (\ref{52})
une suite exacte~:
\begin{equation}
  \label{53}
  0\to\H^0(\maF_s\tens\maG_s)\to\H^1(\maA_s\tens\maG_s)\stackrel{a_s}{\to}\H^1(\maB_s\tens\maG_s)\to\H^1(\maF_s\tens\maG_s)\to 0
\end{equation}
et
$h^0(\maA_s\tens\maG_s)=h^0(\maB_s\tens\maG_s)=h^2(\maA_s\tens\maG_s)=h^2(\maB_s\tens\maG_s)=0.$
Ces annulations et le fait que les familles $\maA\tens\maG$ et $\maB\tens\maG$ sont plates
sur $S$ nous assurent que
$R^1pr_{1*}(\maA\tens\maG)$ est localement libre de rang
$h^1(\maA_s\tens\maG_s)=-\chi(\maA_s\tens\maG_s)$. Le m\^eme argument montre que $R^1pr_{1*}(\maB\tens\maG)$ est localement libre de rang
$-\chi(\maB_s\tens\maG_s)$. La suite (\ref{53}) implique
$$\chi(\maA_s\tens\maG_s)-\chi(\maB_s\tens\maG_s)=\chi(\maF_s\tens\maG_s).$$
Le lemme \ref{lemme1} c) implique
$[\maF_s\tens\maG_s]=[\maF_s]\cdot[\maG_s]$ dans $\kp$. Comme $c$ et
$u$ sont orthogonales, on a $\chi(\maF_s\tens\maG_s)=0$. Alors les \fx
$R^1pr_{1*}(\maA\tens\maG)$ et $R^1pr_{1*}(\maB\tens\maG)$ ont m\^eme
rang. $\Box$

\vs En utilisant la suite exacte (\ref{52}) on d\'efinit un \fib
$\D_S$ par
$$\D_S:=[\det pr_{1!}(\maF\cdot\maG)]^{(-1)}=\det
R^1pr_{1*}(\maB\tens\maG)\tens [\det R^1pr_{1*}(\maA\tens\maG)]^{(-1)}.$$
L'application $a$ fournit une \sct $\sigma_S$ de ce \fib \inve sur
$S$. Ni le \fib $\D_S$, ni la \sct
$\sigma_S$, \`a une fonction inversible  pr\`es, ne d\'ependent de la r\'esolution
choisie. La suite exacte (\ref{53}) montre que la section $\sigma_S$
s'annule exactement aux points $s\in S$ o\`u l'application lin\'eaire $a|_s$ n'est
pas inversible. Ces points sont ceux o\`u $
h^0(\maF_s\tens\maG_s)=h^1(\maF_s\tens\maG_s)\ne 0$. 

On a vu que la suite (\ref{52}) \'etait compatible avec les changements
de base
$S'\to S$. Il r\'esulte que le \fib \inve $\D_S$ et la \sct $\sigma_S$
le sont aussi. Soit $\phi:S\to\emceu$ le \mor modulaire associ\'e aux
familles $\maF, \maG$.

\begin{prop}
  \label{lemme3}
Il existe un \fib \inve $\decu$ sur $\emceu$ et une section
$\sigmaceu$ bien d\'etermin\'ee \`a une constante multiplicative pr\`es, qui v\'erifie: pour toute vari\'et\'e
\alg $S$ et pour toutes familles plates $\maF$, $\maG$
de \fx \sss de classes $c$ \resp $u$ dans $\kp$, param\'etr\'ees par $S$,
on a
\begin{eqnarray*}
  \dees&=&\phi^*(\decu) {\rm\ \ et \ \ }\\
 \sigmaes&=&\phi^*(\sigmaceu) \mbox{\rm \ \ \`a une fonction inversible pr\`es.\ \ }
\end{eqnarray*}
\end{prop}

\par{\bf Preuve~:}

La preuve est classique. On commence par les deux lemmes suivants.

\begin{lemme}[\cite{Simp}, \cite{LeP-Jussieu}, \cite{LeP-Nice}]
\label{lemme:omegace}
Il existe une vari\'et\'e lisse $\Omegace$, une famille plate $\efce$ de
\fx \sss de classe $c$ sur $\pp$ param\'etr\'ee par $\Omegace$ et un
groupe r\'eductif $\gece$ qui agit sur $\Omegace$ tels que $\emce$ soit
un bon quotient de $\Omegace$ sous l'action de $\gece$.
  
\end{lemme}

La description de $\Omegace$ est la suivante. On introduit  pour $m$
entier assez grand la \cep $P(m)$ de $c(m)$, et la somme directe $B$
de $N=P(m)$ exemplaires du \fib \inve $\O_{\pp}(-m)$. On consid\`ere le
sch\'ema de Hilbert $\Hilb^c(B)$ des \fx coh\'erents $F$ de classe de \Gro
$c$ quotients de $B$. Le groupe $\gece:={\rm Aut}(B)$ op\`ere de mani\`ere
naturelle sur $\Hilb^c(B)$. On note  $\Omegace:=\Omega^{ss}$ l'ouvert des points \sss
pour l'action de $\gece$. Ces points
repr\'esentent les \fx \sss $F$ pour lesquels le \mor naturel
$\H^0(B(m))\to\H^0(F(m))$ est un isomorphisme. L'ouvert $\Omega^{ss}$ est
invariant par l'action du groupe r\'eductif ${\rm Aut}(B)$. C'est un ouvert  lisse et
 sur $\Omega^{ss}\times\pp$ on dispose d'un \fs quotient universel
 $\maF$, lequel est aussi muni d'une action de $\gece$. L'ouvert
 $\Omega^{ss}$ admet pour  bon quotient l'espace de modules grossier $\emce$ des \fx
semi-stables de classe $c$.

\begin{lemme}[\cite{LeP2}]
\label{lemme:omegau}
  Il existe une vari\'et\'e lisse $\Omegau$, une famille plate $\gheu$ de
\fx \sss de classe $u$ sur $\pp$ param\'etr\'ee par $\Omegau$ et un
groupe r\'eductif $\geu$ qui agit sur $\Omegau$, tels que $\emu$ soit
un bon quotient de $\Omegau$ sous l'action de $\geu$.
  
\end{lemme}

La description de la vari\'et\'e $\emu$ est analogue \`a celle de la vari\'et\'e
$\emce$. Il existe
un entier  $m$ suffisamment grand pour que tout \fs \ss $G$ sur $\pp$
de classe $u=(0,d,\chi)\in\kp$ v\'erifie: le \fs $G(m)$ est engendr\'e par ses
sections globales et $\H^1(G(m))=0$. On consid\`ere $H$ un espace
vectoriel de dimension $n=dm+\chi$ sur $\comp$. Soit $\Hilb^{u}(B)$ le
sch\'ema de Hilbert-Grothendieck des \fx quotients de $B=H\tens\O(-m)$
de classe $u\in\kp$. Le  groupe $\geu:={\rm GL}(H)$ op\`ere de mani\`ere
naturelle sur $\Hilb^{u}(B)$. On consid\`ere  l'ouvert $\Omegau:=\Omega^{ss}$ des points \sss
pour l'action de $\geu$: ces points correspondent aux \fx quotients de $B$ qui sont
\sss et tels que le \mor d'\'evaluation
$H\to\H^0(F(m))$ soit un isomorphisme. Cet ouvert est lisse et
 $\emu$ est le bon quotient de $\Omega^{ss}$ pour l'action du  groupe $\geu$.

\vs On note $\rho$ la projection $\Omegaceu\to\emceu$. On consid\`ere la
construction pr\'ec\'edente de $\dees$ dans le cas o\`u $S:=\Omegaceu,
\maF:=pr^*_{13}(\efce), \maG:=pr^*_{23}(\gheu)$, o\`u $pr_{13}:\Omegaceu\times\pp\to\Omegace\times\pp$ et
$pr_{23}:\Omegaceu\times\pp\to\Omegau\times\pp$ sont les
projections. On note $\deomega$ le \fibd sur $\Omegaceu$ ainsi obtenu
et $\sigmaom$ sa section canonique, \`a une fonction inversible pr\`es.

De la compatibilit\'e du \fib $\dees$ aux changements de base $S'\to S$,
il r\'esulte une action du groupe $\geceu$ sur le \fib $\deomega$. La
section $\sigmaom$ est \'equivariante. Dans le
lemme qui suit on v\'erifie que la condition de descente est satisfaite
pour $\deomega$.

\begin{lemme}
Soit $\esceesu$ un point de $\Omegaceu$ tel que l'orbite $\gece\cdot
s_c\times\geu\cdot s_u$ soit ferm\'ee dans $\Omegaceu$. Alors le
stabilisateur $\geesceu$ du point $\esceesu$ agit trivialement sur la
fibre $\deesceu$ du \fib $\deomega$ au point $\esceesu$.  
\end{lemme}

\par{\bf Preuve~:}

Les points $\esceesu$ d'orbite ferm\'ee sont les points pour lesquels
l'orbite $\gece\cdot
s_c$ est ferm\'ee dans $\Omegace$ et l'orbite $\geu\cdot s_u$ est ferm\'ee
dans $\Omegau$. D'apr\`es  le lemme 4.2 de \cite{D-N}, le \fs $\efesce$
sur $\pp$ est somme directe de \fx stables de m\^eme polyn\^ome de
Hilbert r\'eduit, et de m\^eme pour $\gheesu$. \'Ecrivons~:
\begin{eqnarray*}
  \efesce&=&F_1^{m_1}\oplus\cdots\oplus F_k^{m_k}\\
 \gheesu&=&G_1^{n_1}\oplus\cdots\oplus G_l^{n_l}
\end{eqnarray*}
pour des \fx stables $F_i, G_j$ diff\'erents deux par deux. Le
stabilisateur du point $\esceesu$ est 
$$\geesceu=\gl(m_1)\times\cdots\times\gl(m_k)\times\gl(n_1)\times\cdots\times\gl(n_l).$$
D'apr\`es la d\'efinition de $\dees$ on a
$$\deesceu=[\det\H^0(\pp,\efesce\tens\gheesu)]^{-1}\tens\det\H^1(\pp,\efesce\tens\gheesu).$$
On a \'egalement
$$\H^q(\pp,\efesce\tens\gheesu)=\oplus_{i=1}^{k}\oplus_{j=1}^{l}\H^q(\pp,F_i^{m_i}\tens
G_j^{n_j})$$
pour $q=0,1$. Par cons\'equent, l'\'el\'ement
$(g_1,\ldots,g_k,h_1,\ldots,h_l)$ appartenant \`a $\geesceu$, agit sur
$\deesceu$ par multiplication avec
$$\prod_{i=1}^{k}\prod_{j=1}^{l}(\det g_i\cdot\det
h_j)^{-h^0(\pp,F_i\tens G_j)+h^1(\pp,F_i\tens G_j)}.$$
Les \fx $F_i, G_j$ ont le m\^eme polyn\^ome de Hilbert r\'eduit que
$\efesce$ \resp $\gheesu$. Puisque les classes $c$ et $u$ sont
orthogonales, on obtient $h^0(\pp,F_i\tens G_j)=h^1(\pp,F_i\tens
G_j)$. Donc le stabilisateur du point $\esceesu$ agit trivialement sur
$\deesceu$.$\Box$

\vs Le lemme de Kempf (lemme de descente, th. 2.3, \cite{D-N}) implique
l'existence d'un unique \fib \inve $\decu$ sur $\emceu$ qui satisfait
$\rho^*(\decu)=\deomega$. Mais $\emceu$ est un bon quotient de
$\Omegaceu$ par l'action du groupe $\geceu$, donc 
$$\H^0(\emceu,\decu)=\H^0(\Omegaceu,\deomega)^{\geceu}.$$
On peut choisir une r\'esolution (\ref{50}) $\gece$-\'equivariante. La
suite exacte (\ref{52}) sera alors $\geceu$-\'equivariante. Il
r\'esulte que la section $\sigmaom$, bien d\'etermin\'ee  \`a
fonction inversible pr\`es,
est $\geceu$-\'equivariante, et donc bien d\'etermin\'ee  \`a une
constante multiplicative pr\`es. 
On d\'eduit l'existence d'une section
$\sigmaceu\in\H^0(\emceu,\decu)$, bien d\'etermin\'ee  \`a une
constante multiplicative pr\`es.

Le fait que le \fib $\decu$ et la section $\sigmaceu$ satisfont la
propri\'et\'e d'universalit\'e de l'\'enonc\'e r\'esulte d'un argument classique
(\cite{LeP-Bergen}, \S 2.13). Cela termine la preuve de la proposition \ref{lemme3}.$\Box$
\begin{lemme}
  Le \fib \inve $\decu$ sur $\emceu$ est isomorphe au \fib \inve $\deuce$.
\end{lemme}

\par{\bf Preuve~:}

Prouvons que la restriction $\deceesu$ du \fib $\decu$ \`a $\emceesu$
est isomorphe au \fib $\deu$ pour tout $s\in\emu$. Soit $G$ un \fs sur
$\pp$ dans la classe du point $s_u$. Soit $S$ une vari\'et\'e et $\maF$
une famille plate de \fx \sss de classe $c\in\kp$, param\'etr\'ee par
$S$. On note $\varphi:S\to\emce$ le \mor modulaire associ\'e \`a $\maF$ et
$\phi=\varphi\times s_u:S\to\emceu$ le \mor modulaire associ\'e \`a $\maF, G$.
D'apr\`es la propri\'et\'e d'universalit\'e on a
$$\varphi^*\deceesu=\phi^*\decu=\det pr_{1!}(\maF\tens pr_2^*(G))^{-1}.$$
Alors il y a un \iso $\varphi^*\deceesu\simeq\varphi^*\decu$ pour tout
couple $(S,\maF)$. On d\'eduit que $\deceesu$ et $\deu$ sont
isomorphes. Nous prouvons de la m\^eme mani\`ere que la restriction de
$\decu$ \`a $\esceemu$ est isomorphe au \fib $\dece$ sur $\emce$. Les
vari\'et\'es $\emu$ et $\emce$ sont projectives et int\`egres. Le lemme
suivant s'applique.

\begin{lemme}[\cite{Mum}, p. 23, \cite{Milne}, \S 5, th. 5.1 et cor. 5.2]
  \label{lemme3.5}
Soit $M$ une va\-ri\-\'e\-t\'e \alg in\-t\`egre, $N$ une  vari\'et\'e \alg projective,
et $\maL$ un \fs \inve sur
$M\times N$. On suppose que la classe d'isomorphisme de la restriction de $\maL$ \`a la fibre
$\{m\}\times N$ est la m\^eme pour chaque $m\in M$. Alors 
$\maL$ s'\'ecrit $L_1\boxtimes L_2$ pour deux \fx $L_1\in\Pic(M)$,
$L_2\in\Pic(N)$. 

\end{lemme}
\vs Du lemme et du calcul des restrictions du \fib $\decu$ aux
 fibres on trouve $\decu\simeq\deuce.$ $\Box$

\vs Cela prouve le point (i)  du \th \ref{pro:cadru}. Les points
(ii) et (iii) sont \'evidents. Le point (iv) r\'esulte du lemme g\'eom\'etrique \'evident~:
\begin{lemme}
  \label{lemme4}
Soient $M$ et $N$ deux vari\'et\'es projectives, $\D$ et $\maE$ des
\fibs \inves sur $M$ \resp $N$, et  une \sct
$\sigma\in\H^0(M\times N,\D\boxtimes\maE)$. On note
$\sigma_m\in\Gamma(N,\maE)$ la restriction de $\sigma$ \`a $\{m\}\times
N$. On suppose que $\sigma_m$ n'est pas identiquement nulle pour tout
$m\in M$. Alors~:

i) La \sct $\sigma$ produit un \mor
$D_{M,N}:\H^0(N,\maE)^*\to\H^0(M,\D)$.

ii) La \sct $\sigma$ produit une application $\Phi:
M\to\proj\H^0(N,\maE)$ d\'efinie par $m\mapsto[\sigma_m]$ qui v\'erifie
$\Phi^*\O(1)=\D$.

iii) L'\app $\Phi$ induit sur les sections globales une \app
$\Phi^*:\H^0(N,\maE)^*\to\H^0(M,\D)$. Alors $\Phi^*=D_{M,N}$.

iv) Si l'image du \mor $\Phi$ n'est pas contenue dans un hyperplan
(c'est-\`a-dire que les $\sigma_m$ engendrent $\H^0(N,\maE)$),
alors $D_{M,N}$ est injectif.

\end{lemme}

On applique le lemme pour $M=\emce, N=\emu, \D=\deu, \maE=\dece,
\sigma=\sigmaceu$, lorsque $\sigma_F$ n'est pas identiquement
nulle pour tout $[F]\in\emce$. $\Box$

\begin{rems}
  
{\rm 1) D'apr\`es le point i) du \th \ref{pro:cadru}, le fait que $\sigma_F$ n'est pas identiquement
nulle \'equivaut \`a: il existe
$G\in\emu$ tel que $h^0(F\tens G)=h^1(F\tens G)=0$. Dans les
situations consid\'er\'ees cela sera vrai  pour tout $F\in\emce$. 
Le Potier a montr\'e  cette affirmation
 si  $2c_1^c=0 {\rm \ mod\,} r^c$, en utilisant
l'existence d'une droite, ou bien d'une conique qui n'est pas de saut
pour un \fs stable g\'en\'erique $F$ de classe $c$. 

2) Dans  la proposition
3.3 de \cite{Kyoto}, pour une classe $c$ donn\'ee telle que l'espace
$M_c$ est non-vide, on montre l'existence de $u\in c^{\bot}$ de
dimension $1$ telle que $\sigmaceu\ne 0$. La d\'emonstration repose sur
un th\'eor\`eme de Flenner, qui donne le comportement de la semi-stabilit\'e
par restriction aux courbes de degr\'e \'elev\'e, et la version effective
d'un r\'esultat de Faltings.}
\end{rems}

\section{Pr\'eliminaires: le \fs canonique sur $\emd$ et sur $\emce$}

\label{premiers}

Le but de cette section est de donner des descriptions explicites pour
les \fx canoniques sur $\emd$ et sur $\emce$. Ces r\'esultats, et le
\th d'annulation de Kawamata-Viehweg, seront appliqu\'es pour
d\'eduire des r\'esultats d'annulation (thm. \ref{the:pou}, prop. \ref{propo2}).

Nous allons \'etudier $\emu=\emd$, pour $u=d\gotu=d(0,1,0)$ et $d\ge 1$. On d\'efinit
$C_d=\proj\H^0(\pp,\O(d))$, l'espace des courbes de degr\'e $d$ dans
$\pp$. Pour un \fs 
$G\in \emu$, il existe une pr\'esentation
\begin{equation}
  \label{58}
 0\to Q\stackrel{a}{\to}R\to G\to 0 
\end{equation}
et $\det a\in\H^0(\pp,\det R\tens\det Q^{-1})=\H^0(\pp,\O(d))$
s'appelle l'\'equation du support sch\'ematique de $G$. Cela ne d\'epend
pas, \`a une constante pr\`es, du choix de $Q,R$. On d\'efinit ainsi une
\app 
\begin{equation}
  \label{59}
\pi:\emu\to C_d.  
\end{equation}

\subsection{ Le \fs canonique sur $\emd$}
\label{omegauri}
La description de la vari\'et\'e $\emd$ comme quotient d'un ouvert lisse
$\Omega^{ss}$ d'un sch\'ema de Hilbert-Grothendieck par l'action d'un
groupe r\'eductif $\sl(H)$ a \'et\'e donn\'ee au paragraphe \ref{morphisme},
lemme \ref{lemme:omegau}. 
\begin{prop}
  \label{omegau}
Le \fs dualisant de la vari\'et\'e $\emu$,
 $\omega=\omegau$, est inversible et isomorphe \`a $\pi^*\O(-3d)$.
\end{prop}

\par{\bf Preuve~:}

On peut supposer $d\ge 3$, puisque pour $d=1,2$, $\pi$ est un
isomorphisme (voir prop. \ref{pro:d12}) et l'\'egalit\'e est v\'erifi\'ee. 
Il r\'esulte du
th\'eor\`eme de Boutot, \cite{Bou}, que $\emu$ est une vari\'et\'e \`a
singularit\'es rationnelles. En particulier c'est une vari\'et\'e normale et
de
Cohen-Macaulay. 
Soient $\em_u^s$ l'ouvert de $\emu$ des classes de \fx stables, $j$
l'inclusion canonique $j:\em_u^s\to\emu$ et $Y=C_{\em_u^s}$ le
compl\'ementaire de $\em_u^s$  dans $\emu$.  Il est d\'emontr\'e
dans \cite{LeP2}, prop 3.4, que $\codim Y\ge 2$. Soit $\underp\in Y$ 
le point g\'en\'erique d'une sous-vari\'et\'e irr\'eductible de $Y$.
 Puisque $\emu$ est de
Cohen-Macaulay on obtient $\profo(\omega_{\underp})\ge 2$. Or on a l'\'enonc\'e
suivant
(\cite{grot})~:
\begin{theor}
  Soit $X$ un sch\'ema et $Y\subset X$ un ferm\'e. Soit $F$ un \fs \alg
  \co sur $X$ dont le support est $X$. Soit $n$ un entier. Les
  conditions suivantes sont \'equivalentes~:

i) Pour tout point $x\in Y$, on a
$$\profo(F_x)\ge n.$$

ii) Pour $i<n$
$$\has^i_Y(F)=0.$$
\end{theor}
Cet \'enonc\'e entra\^\i ne que $\has^i_Y(\emu,\omega)=0$ pour $i=0,1$. De la
suite longue de cohomologie \`a support on d\'eduit que
$\omega=j_*(j^*(\omega))$. Le m\^eme argument appliqu\'e au \fs inversible
$\pi^*\O(-3d)$ implique $\pi^*\O(-3d)=j_*(j^*(\pi^*\O(-3d)))$, donc il
suffit de d\'emontrer l'\iso souhait\'e sur l'ouvert $\em_u^s$.

On note $\Omega^{s}$ la pr\'eimage de $\em_u^s$ par le \mor $\rho:\Omega^{ss}\to
\emu$.

\begin{lemme}
  \label{lem:rhostar}
L'application $\rho^*:\Pic(\em_u^s)\to\Pic(\Omega^s)$ est injective.

\end{lemme}

\par{\bf Preuve~:} 

L'action de $\sl(H)$ sur $\Omega^{ss}$ se factorise
\`a travers une action propre et libre du groupe $G={\rm PSL}(H)$, et
$\em_u^s$ est le quotient de cette action (cf. \cite{LeP2}, lemme
2.4). En appliquant le lemme de descente de Kempf (th. 2.3, p. 63, et
remarque p. 66, \cite{D-N}), on obtient un isomorphisme
$\Pic(\em_u^s)=\Pic^G(\Omega^s)$, o\`u $\Pic^G$ d\'esigne le groupe des
fibr\'es inversibles munis d'une action de $G$. On a la suite exacte
(\cite{LeP-Durham}, \S 3.3)~:
$$ 0\to \H^1(G,\O^*(\Omega^s))\to\Pic^G(\Omega^s)\to\Pic(\Omega^s)$$
o\`u $\H^1(G,\O^*(\Omega^s))$ est  l'espace des \mors crois\'es
$\phi:G\times\Omega^s\to \comp^*$ , c'est-\`a-dire qui v\'erifient
$$\phi(gg',x)=\phi(g,g'x)\cdot\phi(g',x).$$

 Mais  $G={\rm PGL }(H)$ et les seules fonctions r\'eguli\`eres inversibles sur
${\rm GL }(H)$  sont les caract\`eres de ${\rm GL }(H)$, donc les seules
fonctions r\'eguli\`eres inversibles sur ${\rm PGL }(H)$ sont les
constantes.  Pour un \mor crois\'e $\phi:G\times\Omega^s\to\comp^*$ on a
$\phi(e,x)=1$ et la fonction $\phi_x:G\to\comp^*$ est r\'eguli\`ere
inversible. Cela implique que tout morphisme crois\'e est constant \'egal
\`a $1$, et la conclusion. $\Box$

\vs Par cons\'equent, il suffit de montrer 
\begin{equation}
  \label{24}
  \rho^*\omega= \rho^*\pi^*\O(-3d)
\end{equation}
dans $\Pic(\Omega^s)$. Il suffit encore de le prouver dans
$\Pic(\Omega^s)\tens\qiu$, puisque $\Pic(\emu)=\Pic(\em_u^s)$ est sans
torsion (cf. th. 3.5, \cite{LeP2}). On d\'emontre cette affirmation en
appliquant la formule de Riemann-Roch-Grothendieck. Si $\tang_{\em_u^s}$ est
le \fib tangent \`a $\em_u^s$ et $\maG$ est la famille universelle de \fx
stables
de dimension $1$ param\'etr\'ee par $\Omega^s$, on a
$\rho^*\tang_{\em_u^s}=\uExt^1_{pr_1}(\maG,\maG).$ Puisque pour chaque
$s\in \Omega^s$, $\maG_s$ est un \fs stable, on a
$\Hom(\maG_s,\maG_s)=0$ pour tout $s$, donc
$\uHom(\maG,\maG)=\O$. Alors
$\uExt^0_{pr_1}(\maG,\maG)=pr_{1*}\uHom(\maG,\maG)=pr_{1*}\O=\O_{\Omega^s}$.
Chaque \fs $\maG_s$ est de dimension $1$, donc on a aussi
$\uExt^i_{pr_1}(\maG,\maG)=0$ pour $i\ge 2$. Alors 
$$\det\uExt^{\bullet}_{pr_1}(\maG,\maG)=[\det\uExt^1_{pr_1}(\maG,\maG)]^{-1}.$$
On peut calculer $\rho^*\omega=\rho^*(\det
\tang_{\em_u^s})^{-1}=(\det\rho^*\tang_{\em_u^s})^{-1}=\det\uExt^{\bullet}_{pr_1}(\maG,\maG)$
dans $\Pic(\Omega^s)\tens\qiu$ avec la formule de
Riemann-Roch-Grothendieck~:
$$\rho^*\omega=pr_{1*}([(ch\maG)^*(ch\maG)Td(\pp)]_3).$$
On a 
\begin{equation}
  \label{25}
[(ch\maG)^*(ch\maG)Td(\pp)]_3=-(ch_1\maG)^2Td_1(\pp)=\frac{1}{2}c_1^2(\maG)c_1(\omega_{\pp}).
\end{equation}

Calculons $c_1(\maG)\in\Pic(\Omega^s\times\pp)$. Par le corollaire
\ref{corcor} on a 
$c_1(\maG)=\rho^*\pi^*\O(1)\boxtimes\O(d)$. On note
$\gothas=\rho^*\pi^*\O(1)$ et $h=pr_{2}^*\O(1)$ dans
$\Pic(\Omega^s\times\pp)$. En notation additive on trouve
$c_1(\maG)=\gothas+dh$. On obtient dans (\ref{25})~:
$$\rho^*\omega=pr_{1*}([\frac{1}{2}(\gothas+dh)^2(-3h)])=-3d\gothas=\rho^*\pi^*\O(-3d)$$
par la formule de projection.
L'\'egalit\'e (\ref{24}) est prouv\'ee, donc la proposition.$\Box$ \hfill $\Box$

\subsection{Le \fs canonique et le \fibd de Donaldson sur $\emce$}
\label{canodona}
\vs Le long de ce paragraphe  $c$ d\'esignera une classe g\'en\'erale $c=(r,c_1,\chi)$  dans $\kp$ telle que
  $r>0$ et $\emce$ soit non-vide. Le \fib $\D=\deunu$ sera le \fibd de
  Donaldson associ\'e \`a la
  classe orthogonale \`a $c$,
  $\gotu=(0,\frac{r}{\delta},-\frac{c_1}{\delta})$ o\`u $\delta={\rm
  pgcd (r,c_1)}$. 

\begin{prop}
  \label{omegace}  Le \fs dualisant
  $\omegace$ est inversible et on a $\omegace\simeq\D^{\otimes -3\delta}$ dans $\Pic(\emce)$.
\end{prop}

\par{\bf Preuve~:}

 Soit  $\em_c^s$ l'ouvert des classes repr\'esentant des \fx
 stables. Il y a deux cas \`a consid\'erer: quand la codimension du
 compl\'ementaire $C_{\em_c^s}$ de  $\em_c^s$ dans $\emce$ est $\ge
2$ et quand ce ferm\'e est une hypersurface.

Dans le premier cas la d\'emonstration est identique \`a celle de la
proposition \ref{omegau}. 
On se ram\`ene \`a prouver l'\iso sur $\em_c^s$. Comme le
groupe 
$\Pic(\emce)=\Pic(\em_c^s)$ est sans torsion, il suffit de prouver
l'\iso dans $\Pic(\Omega^s)\tens\qiu$, o\`u $\Omega^s$ est l'image
r\'eciproque de $\em_c^s$ dans $\Omega^{ss}$. Les calculs pour les
premi\`eres classes de Chern des fibr\'es $\rho^*(\omega_{\emce})$ et
$\rho^*(\D^{\tens -3\delta})$ dans $\Pic(\Omega^s)\tens\qiu$ ont \'et\'e
faits par O'Grady  (\cite{O'Gra}) en utilisant la formule de
 Riemann-Roch-Grothendieck.
Puisque ces classes co\"\i ncident, le r\'esultat d\'ecoule. 

Une analyse du second cas a \'et\'e faite par Dr\'ezet (\cite{Drez}). Il est
prouv\'e que~: la classe $c$ est divisible par $2$, l'espace de modules
$\emcepedoi$ s'identifie \`a $\pp$, le \fibd $\deunu$ sur
$\emcepedoi$ s'identifie \`a $\O_{\pp}(1)$ et l'espace de modules $\emce$
s'identifie \`a $\proj_5$. L'hypersurface des points strictement
semi-stables est l'image du \mor
 $$\sym^2(\emcepedoi)\to \emce$$
qui associe  aux classes $[E],[F]$ la classe $[E]\oplus[F]$ dans
$\emce$. Lorsque la classe $[F]$ est fix\'ee, le \mor
$$\phi_F:\emcepedoi=\pp\to\emce=\proj_5$$
qui associe \`a  $[E]$ la classe $[E]\oplus[F]$ est lin\'eaire.

Cela suffit pour terminer la preuve de la proposition. En effet, le
\fib canonique sur $\proj_5$ est $\omega_{\proj_5}=\O_{\proj_5}(-6)$
et il suffit de prouver que $\deunu$ s'identifie \`a
$\O_{\proj_5}(1)$. On prouve facilement que
$\phi^*_F(\deunu)=\deunu$. Puisque $\phi_F$ est lin\'eaire on a aussi
$\phi^*_F(\O_{\proj_5}(1))=\O_{\pp}(1)$. On avait vu que l'\iso
$\deunu\simeq\O_{\pp}(1)$ \'etait satisfait sur
$\emcepedoi=\pp$. Puisque l'application $\phi^*_F$ est bijective on
obtient que $\deunu=\O(1)$ dans $\Pic(\emce)$.
$\Box$

\begin{rem}
  {\rm Le calcul du faisceau dualisant sur $\emce$ a d\'ej\`a \'et\'e fait par
  Dr\'ezet \`a partir des monades, sans qu'il ait reconnu le r\^ole du \fibd de Donaldson
  (\cite{Dre}, th. F), et par O'Grady (\cite{O'Gra}), sur l'ouvert des
  points stables. Le calcul fait dans la proposition \ref{omegau} pour le \fs dualisant de $\emd$
  est calqu\'e sur ce dernier.}
\end{rem}

\begin{prop}
\label{nefetbig}
Le fibr\'e $\D$ est nef et big.
\end{prop}

\par{\bf Preuve:}
 
D'apr\`es (\cite{LeP-Bergen}, \cite{Li}) il existe un entier $k$ satisfaisant
  aux conditions suivantes:

-le fibr\'e $\D^{\otimes k}$  est engendr\'e par ses sections;

-consid\'erons le \mor associ\'e
$$\phi_k:\emce\to\proj_{\bullet}\H^0(\emce, \D^{\otimes k})$$
dans l'espace projectif des hyperplans de $\H^0(\emce, \D^{\otimes
  k})$. La restriction de $\phi_k$ \`a l'ouvert des fibr\'es $\mu$-stables
est \`a fibres finies.

De la premi\`ere condition on d\'eduit que $\D$ est nef, de la seconde  que le nombre $\int_{\emce}c_1(\D)^{\dim \em_c}$ est
strictement positif, c'est-\`a-dire que $\D$ est big.$\Box$

\vs Par un \th de Boutot (\cite{Bou}), l'espace de modules $\em_c$ est \`a
singularit\'es rationnelles. Par ailleurs le \th de Kawamata-Viehweg est
valable sur les vari\'et\'es  \`a
singularit\'es rationnelles (\cite{E-V}). Nous obtenons d'apr\`es les
propositions \ref{omegace} et \ref{nefetbig} le
th\'eor\`eme~:    
\begin{theor}
\label{the:pou}
Pour $q>0$ et $k>-3\delta$, on a $\H^q(\em_c,\D^{\otimes k})=0$.
\end{theor}

\section{Les espaces de modules $\emd$}
\label{lesespaces}

On consid\`ere la classe d'un point $\eta^2=(0,0,1)\in\kp$. C'est une
classe orthogonale \`a $u=d\gotu$.
Dans la suite, nous allons \'etudier les \fibs $\dece$ sur $\emu$, pour
$\langle c,u \rangle=0$. On peut \'ecrire $c$ sous la forme
$r(c)[\O]-n\eta^2\in\kp$. Par additivit\'e, il suffit d'\'etudier
$\D_{\eta^2}$ et $\D_{\O}$.

\begin{prop}
  \label{pro:pistar}
On a $\pi^*\O(1)=\D_{\eta^2}^{-1}$.
\end{prop}

\par{\bf Preuve~:} 

On note $\Xi$ l'hypersurface universelle dans
$C_d\times \pp$
param\'etr\'ee par $C_d$.  La proposition 2.8 de \cite{LeP-Durham} affirme
que $\D_{\eta^2}^{-1}=\pi^*(\lambda_{\O_{\Xi}}(\eta^2))$ o\`u
$$\lambda_{\O_{\Xi}}(\eta^2)=\det pr_{1!}(\O_{\Xi}\cdot
pr_2^*(\eta^2)).$$
En partant de
la r\'esolution de $\O_{\Xi}$ sur $C_d\times \pp$~:
$$0\to \O(-1,-d)\to\O\to \O_{\Xi}\to 0$$
on obtient 
$\lambda_{\O_{\Xi}}(\eta^2)=\O(1)$ sur $C_d$. $\Box$ 

\vs Soit $\maG$ une famille de \fx $G$ de dimension $1$ param\'etr\'ee par une
vari\'et\'e \alg int\`egre $S$, et 
$$\phi: S\to\emd$$
 le \mor modulaire
associ\'e. 
\begin{cor}
  \label{corcor}
On a $\det\maG=(\phi\circ\pi)^*(\O(1))\boxtimes\O(d)$ sur $S\times \pp$.
\end{cor}

\par{\bf Preuve du corollaire~:}

On note $(\det\maG)_s$ la restriction de $\det\maG$ \`a $\{s\}\times
\pp$ et $(\det\maG)_x$ la restriction de $\det\maG$ \`a
$S\times\{x\}$. D'apr\`es le lemme \ref{lemme3.5} il suffit de prouver
que $(\det\maG)_s=\O(d)$ sur $\pp$ et que
$(\det\maG)_x=(\phi\circ\pi)^*(\O(1))$ sur $S$. On a \'evidemment
$(\det\maG)_s=\O(d)$ car $\maG$ est une famille de \fx de classe
$c_1=d$ sur $\pp$ param\'etr\'ee par $S$. Consid\'erons la r\'esolution
(\ref{49}) pour $\maG$. On note $\maQ_x$, $\maR_x$, $\maG_x$ les
restrictions de $\maQ$, $\maR$, et \resp $\maG$ \`a $S\times\{x\}$. On a
une suite exacte~:
$$0\to\uTor^{S\times\pp}_1(\maG,pr_2^*\O_x)\to\maQ_x\to\maR_x\to\maG_x\to
0.$$
Chacun de ces \fx a un support inclus dans  $S\times\{x\}$, donc
leur image directe sup\'erieure $R^ipr_{1*}$ est nulle pour $i>0$. Il
r\'esulte que 
$$pr_{1!}(\maG\cdot
pr_2^*\O_x)=[pr_{1*}(\maG_x)]-[pr_{1*}(\uTor^{S\times\pp}_1(\maG,pr_2^*\O_x))]=[\maR_x]-[\maQ_x]$$
o\`u on a identifi\'e $\maQ_x$ avec $pr_{1*}(\maQ_x)$ et $\maR_x$ avec
$pr_{1*}(\maR_x)$. Par la d\'efinition de $\D_{\eta^2}^{-1}$ on trouve~:
\begin{equation}
  \label{82}
  \phi^*\D_{\eta^2}^{-1}=\phi^*\D_{[\O_x]}^{-1}=(\det\maQ_x)^{-1}\tens \det\maR_x.
\end{equation}
D'apr\`es la proposition \ref{pro:pistar} il r\'esulte que 
\begin{eqnarray*}
  (\det\maG)_s&=&(\det\maQ_x)^{-1}\tens \det\maR_x\\
&=&\phi^*\D_{\eta^2}^{-1}\\
&=&(\phi\circ\pi)^*(\O(1)).\Box
\end{eqnarray*}

\vs On d\'efinit le \fib \inve $\Theta$ sur $\emu$ comme
$\Theta=\D_{\O}=\det pr_{1!}(u)^{-1}$. Tout ce qu'on utilisera dans la
suite est
contenu dans la proposition suivante, extraite de \cite{LeP2},
chap.2~:

\begin{prop}[\cite{LeP2}]
\label{propro}
  Le \fib $\Theta$ a une section \cano $\theta$, unique \`a constante
  pr\`es, non identiquement nulle, qui s'annule aux points $G$ tels que
  $h^0(\pp,G)=h^1(\pp,G)\ne 0$.
\end{prop}

\subsection{Les cas $d=1,2$}

\begin{prop}[\cite{LeP2}]
\label{pro:d12} Pour $d=1,2$, le \mor $\pi$ est un isomorphisme, et le
  \fib $\Theta$ est trivial. En conclusion, les espaces
  $\proj\H^0(\emd,\Theta^{\otimes r}(n))$ et $\proj\H^0(C_d,\O(n))$ s'identifient.
 
\end{prop}

\par{\bf Preuve~:}

 Le fait que $\pi$ est un \iso dans ce cas est
d\'emontr\'e dans \cite{LeP2}. Plus pr\'ecis\'ement, l'inverse de $\pi$ est
donn\'e de la mani\`ere suivante: pour $d=1$, \`a une droite $l$ on associe
le \fs $\O_l(-1)$, pour $d=2$ et pour une conique $C$ lisse, \`a $C$ on
associe le \fs $\O_C(-a)\simeq\O_{\proj_1}(-1)$ ($a$ point de $C$),
pour $d=2$ et $C$ d\'ecomposable en deux droites $l_1$ et $l_2$, \`a $C$
on associe le \fs $\O_{l_1}(-1)\oplus\O_{l_2}(-1)$. Pour chacun de ces
\fx on a $h^0=h^1=0$. Il r\'esulte que la \sct $\theta$ est partout non
nulle sur $\emd$. Donc le \fib $\Theta$ est trivial. La conclusion
r\'esulte de la proposition \ref{pro:pistar}, du fait que  $\pi$ est un isomorphisme 
et de la trivialit\'e du \fib $\Theta$. $\Box$

\vs On d\'eduit en corollaire~:
\begin{cor}
  Pour $c=2-n\eta^2$ on a
  $$\H^0(\emunu,\dece)=\H^0(\pps,\O_{\pps}(n))=\es^nE^*$$
   et
$$\H^0(\emdoi,\dece)=\H^0(\proj_5,\O_{\proj_5}(n))=\es^n(\es^2E^*).$$
\end{cor}
On rappelle que $E=\H^0(\pp,\O(1))$.

\subsection{Le cas $d=3$}

Soit $\gotu$ la classe $(0,1,0)\in \kp$ et $c=(2,0,c_2=n)$.
 L'objectif de ce paragraphe est de d\'emontrer la~:
\begin{prop}
\label{propropro}
  En tant que $\sl(3)$-\rep l'espace $\H^0(\emtrei,\dece)$ s'identifie
  \`a $\es^n(\es^3E^*)\oplus\es^{n-2}(\es^3E^*)$ ( o\`u $E=\H^0(\pp,\O(1))$).
\end{prop}

Notre r\'ef\'erence principale sera l'article
\cite{LeP2}. 
On consid\`ere la cubique universelle $\maC\subset C_3\times\pp$. La
projection $\maC\to C_3$ induit par la propri\'et\'e universelle de
l'espace de modules $\emtrei$ un \mor $s:C_3\to\emtrei$ qui associe \`a
une cubique $C$ son \fs structural $\O_C$. Le \mor $s$ est une section
de $\pi$~: la r\'esolution
$$0\to\O(-3)\stackrel{\cdot \uC}{\to}\O\to\O_C\to 0$$
d\'emontre que $(\pi\circ s)(C)=C$. Ici $\uC$ d\'esigne l'\'equation de la
cubique $C$.

Prenons $D\in\pp\times\pps$ la vari\'et\'e d'incidence et $p$ et $q$ les
projections~:
$$\diagram
D\dto_{q}\rto^{p}&\pp\\
\pps.
\enddiagram$$

L'\app $G\to q_*(p^*(G(-2)))(-1)$ d\'efinit (cf. \cite{LeP2}) un \mor
$\phi:\emtrei\to\em^*_{(3,0,0)}$ dans l'\espmod de \fx \sss de classe
$(r=3,c_1=0,\chi=0)$ sur $\pps$. On note encore $\D=\deu$ le \fibd sur
$\em^*_{(3,0,0)}$, associ\'e \`a la classe $\gotu=(0,1,0)$. La proposition
suivante est prouv\'ee dans \cite{LeP2}~:
\begin{prop}
\label{pro:led}
i) Le diviseur des z\'eros $\divs$ de la section $\theta$ co\"\i ncide
avec l'image du \mor $s$.

ii) Le \mor $\phi$ est l'\'eclatement d'un point lisse et la fibre
exceptionnelle est $\divs$.

\end{prop}

On obtient que les \mors $\pi,
s$ \'etablissent un \iso entre $\divs$ et $C_3$. On  identifiera $\divs$
et $C_3$
dans la suite. Puisque $\divs$ est le diviseur
exceptionnel, son \fib conormal est $\O_{C_3}(1)$. De la suite exacte
courte
\begin{equation}
  \label{unurond}
  0\to\Theta^{-1}\stackrel{\theta}{\to}\O\to\O_{\divs}\to 0
\end{equation}
il r\'esulte que
\begin{equation}
  \label{doirond}
  \Theta^{-1}|_{\divs}\simeq\O_{C_3}(1).
\end{equation}

\begin{prop}
\label{thetadeunu}
  On a $\Theta(1)=\phi^*\D$ dans $\Pic(\emtrei)$.
\end{prop}

\par{\bf Preuve~:}

Compte-tenu du fait que $\Pic(\emtrei)$ est sans
torsion (cf. th. 3.5, \cite{LeP2}) il suffit de d\'emontrer que
$$\phi^*\D^{\otimes -9}=\Theta^{\otimes -9}(-9).$$
Par la proposition \ref{omegace}, on a $ \D^{\otimes
  -9}=\omega_{\em^*_{(3,0,0)}}$, et par la proposition \ref{omegau},
on a $\O(-9)=\omega_{\emtrei}$. L'\'egalit\'e \`a d\'emontrer devient
$$\phi^*\omega_{\em^*_{(3,0,0)}}=\omega_{\emtrei}\tens \Theta^{\otimes
  -9}.$$
Mais on a vu que $\Theta=\O(\divs)$ et que $\divs$ est le diviseur
exceptionnel. L'\'egalit\'e \`a d\'emontrer est 
$$\omega_{\emtrei}=\phi^*\omega_{\em^*_{(3,0,0)}}\tens\O(9\,\divs)$$
qui est valable chaque fois qu'on \'eclate un point lisse dans une
vari\'et\'e de dimension $10$ (\cite{Har}, ex. II 8.5, p. 188).$\Box$

\vs Puisque $c=2-n\eta^2$ dans $\kp$, on a $\dece\simeq\Theta^{\otimes 2}(n)$ sur
$\emtrei$. On passe au calcul de l'espace
$\H^0(\emtrei,\Theta^{\otimes 2}(n))$.
\begin{prop}
\label{propo2}
  $\H^q(\emtrei,\Theta(n))=0$ pour $q\ge 1$ et $n\ge -8$.
\end{prop}

\par{\bf Preuve~:}

L'\espmod $\emtrei$ est le quotient d'une vari\'et\'e lisse par un
groupe r\'eductif, donc il est \`a singularit\'es rationnelles, par
un r\'esultat de Boutot \cite{Bou}.  Le \th de
Kawamata-Viehweg s'applique (\cite{E-V}). D'apr\`es la proposition
\ref{nefetbig} le \fib $\D$ est nef et big sur $\em^*_{(3,0,0)}$ et
d'apr\`es 
les propositions \ref{pro:led} et \ref{thetadeunu} le \fib $\Theta(1)$
est l'image r\'eciproque de $\D$ par un \'eclatement. On en d\'eduit que le
\fib $\Theta(1)$ est big et nef sur $\emtrei$. Le \fib 
$\O(n-1)=\pi^*(\O(n-1))$ est globalement engendr\'e pour $n\ge 1$ donc
$\Theta(n)$ est big et nef pour $n\ge 1$. La proposition \ref{omegau} fournit le \fs dualisant sur $\emd$~:
$\omega_{\emd}=\pi^*(\O(-3d))$. Alors
$\Theta(n)\tens\omega_{\emd}^{-1}=\Theta(n+9)$ est big et nef pour
$n\ge -8$. Le r\'esultat en d\'ecoule.$\Box$

\vs On tensorise la suite (\ref{unurond}) par $\Theta^{\otimes 2}(n)$. On obtient,
apr\`es l'identification $\divs=C_3$, et en utilisant l'isomorphisme
(\ref{doirond}), la suite exacte courte sur $\emtrei$~:
$$0\to\Theta(n)\to\Theta^{\otimes 2}(n)\to\O_{C_3}(n-2)\to 0.$$
La proposition \ref{propo2} conduit \`a une suite exacte courte sur les
sections globales, pour $n\ge -8$~:
\begin{equation}
  \label{labelulul}
  0\to\H^0(\emtrei,\Theta(n))\to\H^0(\emtrei,\Theta^{\otimes 2}(n))\to\H^0(C_3,\O_{C_3}(n-2))\to 0.
\end{equation}

\begin{prop}
  \label{propo3}
Soit $u=3\gotu$. Alors les \mors~:
$$\H^0(C_3,\O(n))\stackrel{\pi^*}{\to}\H^0(\emu,\O(n))\stackrel{\cdot\theta}{\to}\H^0(\emu,\Theta(n))$$
sont des isomorphismes.
\end{prop}

\par{\bf Preuve~:} 

On consid\`ere l'ouvert $U\subseteq C_3$ des courbes
irr\'eductibles. Son compl\'ementaire $C_U$ est de codimension $\ge 2$. La
proposition r\'esulte des lemmes suivants~:

\begin{lemme}
\label{lemme13}
  Le compl\'ementaire de l'image r\'eciproque $\pi^{-1}(U)$ de l'ouvert
  $U$ par le \mor $\pi$ est de codimension au moins $2$ dans $\emu$. En plus on
  a l'isomorphisme $\pi_*(\O_{\pi^{-1}(U)})=\O_U$.
\end{lemme}

\begin{lemme}
\label{lemme14}
  Le \mor
  $\pi_*\O_{\emu}\stackrel{\cdot\theta}{\to}\pi_*\Theta$
  est un isomorphisme sur $C_3$.
\end{lemme}

Effectivement, le lemme \ref{lemme13} implique~:
\begin{eqnarray*}
  \H^0(\emu,\O(n))&=\H^0(\pi^{-1}(U),\O(n))=\H^0(U,\pi_*(\O(n)))\\
&=\H^0(U,\O(n))=\H^0(C_3,\O(n)).
\end{eqnarray*}
Il r\'esulte du lemme \ref{lemme14} que
$\pi_*\O(n)\stackrel{\cdot\theta}{\to}\pi_*\Theta(n)$ est un
isomorphisme. En prenant les sections globales sur $C_3$ on obtient
que $\H^0(\emu,\O(n)) \stackrel{\cdot\theta}{\to}\H^0(\emu,\Theta(n))$
est un isomorphisme.

\vs \par{\bf Preuve du lemme \ref{lemme13}~:}

Soit $V\subseteq \emu$ l'ouvert des \fx stables et localement libres sur
leur support. Le lemme 3.2 et la proposition 3.4 de \cite{LeP2}
d\'emontrent que $\codim C_V\ge 2$, o\`u $C_V$ d\'esigne le compl\'ementaire
de $V$. La proposition 2.8 de \cite{LeP2}
affirme que le \mor $V\stackrel{\pi}{\to}C_3$ est lisse. Il r\'esulte
que ses fibres au-dessus de $C_U$ sont de dimension $\dim\emu-\dim
C_3=1$. L'inclusion $\pi^{-1}(C_U)\subset
C_V\cup(V\cap\pi^{-1}(C_U))$ entra\^\i ne $\codim\pi^{-1}(C_U)\ge
2$. Le \th 9 de \cite{Alt} nous assure que le \mor projectif
$\pi^{-1}(U)\to U$ est plat et \`a fibres int\`egres. D'o\`u
$\pi_*(\O_{\pi^{-1}(U)})=\O_U$. $\Box$

\vs\par{\bf Preuve du lemme \ref{lemme14}~:}

On d\'eduit de  la suite exacte (\ref{unurond}) et de l'\iso
\ref{doirond} la suite exacte courte~:
$$0\to\O\stackrel{\cdot\theta}{\to}\Theta\to\O_{C_3}(-1)\to 0$$
sur $\emtrei$. On applique le foncteur $\pi_*$. On obtient la suite
exacte~:
\begin{equation}
  \label{label}
  0\to\pi_*(\O)\to\pi_*(\Theta)\stackrel{\cdot\theta}{\to}\O_{C_3}(-1)\stackrel{\delta}{\to}R^1\pi_*(\O)
\end{equation}
sur $C_3$. 

Soit $W\subset C_3$ l'ouvert des cubiques lisses. La fibre $P_C$ du
\mor $\pi$ au-dessus de $C\in W$ s'identifie \`a la jacobienne de $C$,
${\rm Jac}\,(C)$. La restriction du \fib $\Theta$ \`a $P_C$ s'identifie
au \fib $\Theta$ usuel sur ${\rm Jac}\,(C)$. Alors le \mor
$\H^0(P_C,\O)\to\H^0(P_C,\Theta)$ est un \iso entre des espaces de
dimension $1$. 
Puisque la fibration 
$$\pi^{-1}(W)\to W$$
 est plate, par le \th de semi-continuit\'e, le \mor $\pi_*\O\to\pi_*(\Theta)$ est un \iso de \fibs \inves sur $W$.

Il r\'esulte de la suite
exacte (\ref{label}) que le \mor
$\O_{C_3}(-1)\stackrel{\delta}{\to}R^1\pi_*(\O)$ est injectif sur
$W$. Comme $\O_{C_3}(-1)$ est un \fs inversible, le \mor $\delta$ est
injectif partout sur $C_3$. Par cons\'equent le \mor
$\pi_*\O_{\emu}\stackrel{\cdot\theta}{\to}\pi_*\Theta$ est un \iso
partout sur $C_3$. $\Box$ \hfill $\Box$

\vs D'apr\`es la proposition \ref{propo3} et de la suite exacte
(\ref{labelulul}) on obtient le

\begin{cor}
  \label{cor:***}
On a une suite exacte courte
$$0\to\H^0(C_3,\O(n))\stackrel{\theta^2\cdot\pi^*}{\to}\H^0(\emtrei,\Theta^{\otimes
  2}(n))\stackrel{s^*}{\to}\H^0(C_3,\O(n-2))\to
0$$
pour $n\ge -8$.
\end{cor}

D'o\`u la proposition \ref{propropro}. $\Box$

\section{Injectivit\'e du \mor $\deceu$}
\label{inject}

On commence par regarder le cas o\`u la classe $c$ est de rang $1$. On
utilise ensuite un argument de r\'ecurrence pour \'etendre le r\'esultat au
cas qui nous int\'eresse, o\`u $c$ est de rang $2$. Le d\'ebut de la
r\'ecurrence utilise le cas $c=(1,0,c_2=n)$ \'etudi\'e en pr\'ealable.

\subsection{Le cas $c=(1,0,c_2=n), u=d\gotu$}

\begin{prop}
  \label{propo4} Pour $d=1,2,3$ et $n\ge 0$, l'image du \mor
  $\Phi:\emce\to\proj\H^0(\emd,\dece)$ n'est pas contenue dans un hyperplan.
\end{prop}

On a vu que le \mor 
$$\H^0(C_d,\O(n))\stackrel{\pi^*}{\to}\H^0(\emu,\O(n))\stackrel{\cdot\theta}{\to}\H^0(\emu,\Theta(n))$$
\'etait bijectif pour  $d=1,2,3$. Le lemme suivant sera utile~:

\begin{lemme}
  \label{lemme*} Si $F=I_Z$ est l'id\'eal du sous-sch\'ema $Z$ des $n$
  points distincts $a_1,\ldots,a_n$ de $\pp$,  si $[G]\in\emd$, et
  s'il existe un point  $a_k\in\supp G$, alors
$$h^0(F\tens G)=h^1(F\tens G)\ne 0.$$
\end{lemme}

\par{\bf Preuve~:}

Sans restreindre la g\'en\'eralit\'e, on peut supposer
$a_1,\ldots,a_i\in\supp G$, $a_{i+1},\ldots, a_n\not\in\supp G$, pour
un nombre $i\in\{1,\ldots,n\}$. On tensorise par $G$ la suite exacte~:
$$0\to F\to I_{\hat{Z}}\to\oplus_{j=1}^i\O_{a_j}\to 0$$
o\`u $I_{\hat{Z}}$ d\'esigne l'id\'eal du sous-sch\'ema des $n-i$
  points distincts $a_{i+1},\ldots,a_n$  et $\O_{a_j}$ le \fs structural du point $a_j$. En utilisant l'annulation $\uTor_1(I_ {\hat{Z}},G)=0$
  (puisque $I_{\hat{Z}}$ est trivial au voisinage du support de $G$), on obtient une inclusion $0\to
  \uTor_1(G,\O_{a_i})\to F\tens G$. Mais \`a partir de la r\'esolution de
  longueur $1$ de $G$ par des \fx localement libres $A$ et $B$~:
$$0\to A\stackrel{\alpha}{\to} B\to G\to 0$$
on obtient apr\`es tensorisation par $\O_{a_i}$~:
$$0\to\uTor_1(G,\O_{a_i})\to A|_{a_i}\stackrel{\alpha|_{a_i}}{\to}
B|_{a_i}\to G\tens\O_{a_i}\to 0$$
et $\det\alpha$ est l'\'equation du support de $G$ donc
$\det\alpha|_{a_i}=0$ et $ \uTor_1(G,\O_{a_i})\ne 0$. Le \fs
$\uTor_1(G,\O_{a_i})$ a pour support le point $a_i$, donc
$\H^0(\uTor_1(G,\O_{a_i}))\ne 0$, d'o\`u $\H^0(F\tens G)\ne 0$.$\Box$

\vs On introduit quelques notations. Pour $E=\H^0(\pp,\O(1))$, le point
$a_i$ est un \'el\'ement de $\proj(E^*) $. Alors $a_i^d$ est
un \'el\'ement de $\proj(\es^dE^*) $ et il repr\'esente, \`a une
constante pr\`es, un \'el\'ement dans $\H^0(C_d,\O(1))$. C'est l'\'equation
$H_{a_i}$ de l'hyperplan des courbes de degr\'e $d$ qui passent par le
point $a_i\in\pp$. Alors $H_{a_1}\cdot\cdots\cdot
H_{a_n}\in\es^n\H^0(C_d,\O(1))=\H^0(C_d,\O(n))$.
\begin{lemme}
  \label{lemmevii}
Pour $F=I_Z$ comme dans le lemme \ref{lemme*} on a $\sigma_F=\sigmaceu(F)=cst\cdot\theta\cdot\pi^*(H_{a_1}\cdot\cdots\cdot
H_{a_n})$, pour une constante $cst\in\comp^*$.
\end{lemme}

Ce lemme suffit pour d\'emontrer la proposition \ref{propo4}, puisque
$\{H_{a_i}\}_i$ engendrent $\H^0(C_d,\O(1))$ et les produits de
$\{H_{a_i}\}_i$ engendrent $\es^n\H^0(C_d,\O(1))$.

\vs\par{\bf Preuve du lemme \ref{lemmevii}~:}

Le lemme \ref{lemme*} nous dit que $\sigma_F$ s'annule sur tous les
\fx $G$ dont le support contient le point $a_i$. Ces \fx appartiennent
\`a l'ensemble d'\'equation $\pi^*(H_{a_i})=0$. La section
$\alpha=\frac{\sigma_F}{\prod_{i=1}^{n}\pi^*(H_{a_i})}$ est une  \sct
rationnelle du \fib $\Theta$ sur $\emd$. Puisque $\sigma_F$ s'annule
sur $\pi^{-1}(\{H_{a_i}=0\})$, $\alpha$ est une \sct r\'eguli\`ere. La
proposition \ref{propo3} appliqu\'ee pour $n=0$ nous assure que
$\H^0(\emd,\Theta)=\H^0(\emd,\O)=\comp$. Donc $\alpha=cst\cdot\theta.$
$\Box$

\vs On remarque, au passage, que la dualit\'e \'etrange dans le cas 
  $r^c=1, d=1,2,3$ a \'et\'e prouv\'ee. 
\begin{prop} 
  Le \mor de dualit\'e \'etrange est un \iso dans le cas $c=(1,0,1-n)$,
  $u=(0,d,0)$, pour $d=1,2,3$ et pour un entier positif $n$.
\end{prop}

\par{\bf Preuve~:}

 En effet, le premier membre de la dualit\'e
est $\H^0(\emd,\Theta(n))^*$, isomorphe par la proposition
\ref{propo3} \`a 
$$\H^0(C_d,\O(n))^*=\es^n(\es^dE).$$
 Pour identifier le
second membre on utilise le fait que $\em_{(1,0,c_2=n)}$ co\"\i ncide avec le
sch\'ema de Hilbert $\Hilb^n(\pp)$ des sous-sch\'emas finis de longueur
$n$ de $\pp$.

Soit $\es^n(\pp)$ le quotient de la
puissance $n$-i\`eme $\pp^n$ de $\pp$ par le groupe sym\'etrique
${\mathfrak{S}_n}$. On dispose du \mor de Hilbert-Chow  $HC:\Hilb^n(\pp)\to
\es^n(\pp)$ qui associe \`a un sch\'ema fini $Z$ le cycle
$\sum_{x\in\pp}lg\,Z_xx$. On note $\O(1,1,\ldots,1)^{
  {\mathfrak{S}_n}}$ le quotient du \fib $\O(1,1,\ldots,1)$ par
l'action de ${\mathfrak{S}_n}$ .
\begin{lemme}
\label{lemasim}
  Les \fibs \inves $\deunu$ et $HC^*(\O(1,1,\ldots,1)^{
  {\mathfrak{S}_n}})$ sont isomorphes sur $\Hilb^n(\pp)$.
\end{lemme}

\par{\bf Preuve~:}

On consid\`ere le diagramme~:
$$\diagram
&{\overbrace{\pp\times\cdots\times\pp}^{n}}\dto^{{\mathfrak{S}_n}}\\
\Hilb^n(\pp)\dto_{\Phi}\rto^{HC}&\es^n(\pp)\dto^{\Psi}\\
\proj\H^0(\emunu,\O(n))&\proj\H^0(C_1,\O(n))\lto_{\theta\cdot\pi^*}^{\sim}.
\enddiagram
$$
Ici, $C_1=\pps$. Le \mor $\Psi$ est d\'efini par $\Psi(a_1,\ldots,a_n)=[H_{a_1}\cdots
H_{a_n}]$. Le lemme \ref{lemmevii} prouve que
$$\Phi=(\theta\cdot\pi^*)\circ\Psi\circ HC$$
 sur l'ouvert des points
distincts de $\Hilb^n(\pp)$. Puisque cet ouvert est dense, le
diagramme consid\'er\'e est commutatif.

L'image r\'eciproque du fibr\'e  $\O(1)$ de l'espace
 projectif $\proj\H^0(C_1,\O(n))$, par $\Psi\circ{\mathfrak{S}_n}$, est le fibr\'e $\O(1,1,\ldots,1)$ sur
 $\proj_2^n$  (on peut le v\'erifier sur
 chaque composante). On tient compte de l'isomorphisme (cf. \cite{LeP-Durham},\S 3.4)~:
 $$\Pic(\es^n(\pp))\to\Pic(\proj_2^n)^{ {\mathfrak{S}_n}}.$$
  On obtient que l'image r\'eciproque du
 fibr\'e  $\O(1)$ de $\proj\H^0(C_1,\O(n))$ sur $\es^n(\pp)$
 est  $\O(1,1,\ldots,1)^{ {\mathfrak{S}_n}}$. Par la commutativit\'e du
 diagramme on obtient la conclusion.$\Box$

\vs En passant aux puissances tensorielles sup\'erieures, on trouve
$\deunu^{\tens d}=HC^*(\O(d,d,\ldots,d)^{ {\mathfrak{S}_n}})$. Cela
entra\^\i ne que 
  $\H^0(\hil),\deunu^{\tens d})=\es^n(\es^dE)$. L'injectivit\'e de $\deceu$
  a \'et\'e prouv\'ee dans la proposition \ref{propo4}.$\Box$

\subsection{Le cas $c=(2,0,c_2=n), u=d\gotu$}

\begin{prop}
  \label{prop5} Pour $1\le d\le 3$ et pour $n\ge 2$, l'image du \mor
  $\Phi:\emce\to\proj\H^0(\emd,\dece)$ n'est pas contenue dans un hyperplan.
\end{prop}

La proposition r\'esulte des quatre lemmes suivants~:
\begin{lemme}
  \label{lem:1} La proposition est vraie pour $n=2$.
\end{lemme}
\begin{lemme}
  \label{lem:2} L'application
  \begin{eqnarray*}
  \H^0(\emd,\Theta^{\otimes
  2}(n))\tens\H^0(C_d,\O(1))&\to&\H^0(\emd,\Theta^{\otimes 2}(n+1))
  \\
 s\tens t&\mapsto&s\cdot\pi^*t
  \end{eqnarray*}
est surjective.
\end{lemme}

\begin{lemme}
  \label{lem:3} Soit $n\ge 2$ et $\em_c^0\subset\emce$ l'ouvert des
  points stables. Si l'image $\Phi(\emce)$ n'est pas contenue dans un
  hyperplan, alors $\Phi(\emce^0)$ n'est pas contenue dans un hyperplan.
\end{lemme}
\begin{lemme}
  \label{lem:4}
Soit $F\in\em^0_{(2,0,c_2=n)}$, $x\in\pp$ et $a:F\surto\O_x$ un \mor
surjectif. Alors $F'=\Ker a$ est un \fs \ss et on a
$\sigma_{F'}=\sigma_F\cdot\pi^*H_x$ par l'\app $id\cdot\pi^*$ du lemme
\ref{lem:2}. Ici
$$\sigma_{F'}\in\H^0(\emd,\Theta^{\otimes 2}(n+1)),
\ \ \sigma_F\in\H^0(\emd,\Theta^{\otimes 2}(n))\mbox{\rm\ \  et\ \ } H_x\in\H^0(C_d,\O(1)).$$
\end{lemme}

Les lemmes \ref{lem:1}, \ref{lem:2}, \ref{lem:3} et \ref{lem:4} fournissent une
d\'emonstration par r\'ecurrence de la proposition \ref{prop5}. Le lemme 
\ref{lem:3} dit que les $\sigma_F$, pour $F$ stable, engendrent
$\H^0(\emd,\Theta^{\otimes 2}(n))$. Mais les $H_x$ engendrent $\H^0(C_d,\O(1))$
lorsque $x$ varie. Par le lemme \ref{lem:2} on obtient que les
$\sigma_F\cdot\pi^*H_x$ engendrent
$\H^0(\emd,\Theta^{\otimes 2}(n+1))$. Le lemme \ref{lem:4} assure que de tels
\'el\'ements sont de la forme $\sigma_{F'}$, donc des images par $\Phi$ de
$\em_{(2,0,c_2=n+1)}$. $\Box$

\vs\par{\bf Preuve du lemme \ref{lem:2}~:}

Pour $d=1,2$, le diagramme suivant est commutatif~:
$$\diagram
\H^0(C_d,\O(n))\tens\H^0(C_d,\O(1))\ \ \ \dto\rto^{\ \ \ (\theta^2\cdot\pi^*)\tens
  id\ \ \ }&\ \ \ \H^0(\emd,\Theta^{\otimes 2}(n))\tens\H^0(C_d,\O(1))\dto\\
\H^0(C_d,\O(n+1))\rto^{\theta^2\cdot\pi^*}&\H^0(\emd,\Theta^{\otimes 2}(n+1)).
\enddiagram$$

Par la proposition \ref{pro:d12}, les \mors horizontaux sont des
isomorphismes. Le lemme r\'esulte de la surjectivit\'e du \mor vertical
gauche.

Pour $d=3$, on note $H=\H^0(C_d,\O(1))$.
 Le lemme est  une cons\'equence du diagramme analogue, commutatif~:

$$\diagram
0\rto&\H^0(C_d,\O(n))\tens H\dto\rto&\H^0(\emd,\Theta^{\otimes
  2}(n))\tens H\dto\rto& \H^0(C_d,\O(n-2))\tens H\dto\rto&0\\
0\rto&\H^0(C_d,\O(n+1))\rto^{\theta^2\cdot\pi^*}&\H^0(\emd,\Theta^{\otimes 2}(n+1))\rto^{s^*}&\H^0(C_d,\O(n-1))\rto&0.
\enddiagram$$

Les suites horizontales sont exactes d'apr\`es le corollaire
\ref{cor:***}. Les \mors verticaux lat\'eraux sont surjectifs, donc
aussi le \mor vertical central. $\Box$

\vs\par{\bf Preuve du lemme \ref{lem:3}~:}

Ceci est \'evident, puisque $\em_c^0$ est un ouvert dense de
$\emce$. $\Box$

\vs\par{\bf Preuve du lemme \ref{lem:4}~:}

Pour un sous-\fs $F''$ de rang $1$ de $F'$ on a $c_1(F'')\le 0$
puisque $F''$ est aussi un sous-\fs de $F$, et $F$ est stable. Si $c_1(F'')=0$
alors $\chi(F'')<\frac{2-n}{2}-\frac{1}{2}=\frac{2-(n+1)}{2}$. Donc
$F'$ est semi-stable. 

On tensorise par $G$ la suite exacte:
$$0\to F'\to F\to\O_x\to 0.$$
Soit $[G]\in\emd$. En tenant compte du lemme \ref{lemme1}c), on obtient une suite exacte
$$0\to\uTor_1(G,\O_x)\to F'\tens G\to F\tens G\to G\tens \O_x\to 0.$$
Si $x\in \supp G$ alors $\uTor_1(G,\O_x)\ne 0$ (voir la d\'emonstration
du lemme
\ref{lemme*}) et donc $h^0(F'\tens G)\ne 0$. Sinon  $h^0(F'\tens
G)=h^0(F\tens G)$. On d\'eduit que $\sigma_{F'}=cst\cdot\sigma_{F}\cdot
H_x$. $\Box$

\vs\par{\bf Preuve du lemme \ref{lem:1}~:}

Pour $a,b\in \pp$, la classe du \fs $F=I_a\oplus I_b$ appartient \`a
$\em_{(2,0,c_2=2)}$. On commence par prouver que
$\sigma_F=cst\cdot\theta^2\cdot\pi^*(H_a\cdot H_b)$. La \sct $\sigma_F$
s'annule sur $\pi^{-1}(\{H_a=0\})$ et sur $\pi^{-1}(\{H_b=0\})$
d'apr\`es le lemme \ref{lemme*}. Alors la section rationnelle
$\alpha=\frac{\sigma_F}{\pi^*(H_{a}\cdot H_b)}$ de $\Theta^{\otimes 2}$ est
r\'eguli\`ere. Si $d=1,2$ on applique la proposition \ref{pro:d12} pour
avoir un \iso $\H^0(C_d,\O)=\H^O(\emd,\Theta^{\otimes 2})$. Pour $d=3$ cet \iso
d\'ecoule du  corollaire \ref{cor:***} appliqu\'e pour $n=0$. 
Par cons\'equent
$\alpha=cst\cdot\theta^2$. 

Puisque les produits $H_{a}\cdot H_b$ engendrent
$\H^0(C_d,\O(2))$, le lemme \ref{lem:1} est prouv\'e pour $d=1,2$. 

Pour
$d=3$, en regardant la suite exacte du corollaire \ref{cor:***}, on
montre que les $\sigma_F$ engendrent le sous-espace
$\theta^2\cdot\pi^*\H^0(C_3,\O(2))$ de $\H^0(\emd,\Theta^{\otimes 2}(2))$. Pour
montrer le lemme, en tenant compte du fait que $h^0(C_3,\O)=1$, il
suffit de trouver un \fs $F$ pour lequel $s^*\sigma_F\ne 0$. Ceci
\'equivaut \`a trouver une cubique $C$ pour laquelle
$h^0(F|_C)=h^1(F|_C)=0$. Cette condition est satisfaite pour tout \fs
$F$ localement libre stable et toute cubique $C$. Effectivement, si on
reconsid\`ere la r\'esolution de $\O_C$~:
$$0\to\O(-3)\stackrel{\cdot \underline{C}}{\to}\O\to\O_C\to 0,$$
comme $\uTor_1(F,\O_C)=0$, on obtient
$$0\to F(-3)\stackrel{\cdot \underline{C}}{\to} F\to F|_C\to 0.$$
Alors la suite
$$\H^1(F)\to\H^1(F|_C)\to\H^2(F(-3))$$
est exacte. Le nombre de Hodge $h^1(F)=n-2$ est nul pour $n=2$ et par
la dualit\'e de Serre on a $\H^2(F(-3))=\H^0(F^*)=\Hom(F,\O)$. Ce
dernier groupe est nul, d'o\`u $h^1(F|_C)=0$. L'annulation de
$\Hom(F,\O)$ s'obtient ainsi: S'il existe un \mor non nul $m:F\to\O$,
alors $\im m\subset\O$, donc $\im m=I_Z$, pour $Z$ sous-sch\'ema dans
$\pp$. La stabilit\'e de $F$ implique $c_1(I_Z)=0$ et $\chi(I_Z)>0$,
soit que $I_Z=\O$. Alors $\Ker m$ est localement libre de $c_1=0$ et
$\chi=-1$, ce qui est contradictoire. Donc $\Hom(F,\O)=0$.$\Box$

\section{Preuve de la proposition \ref{4.0.7}}
\label{1}
On reprend la d\'emarche et les notations de l'article
\cite{D1}. Dans cet article, on se fixait un entier
positif $l$, qui dans cette application sera toujours \'egal \`a $1$. On introduit la notion de syst\`eme coh\'erent, qui
consiste \`a consid\'erer en m{\^e}me temps que le \fs $F$, un sous-espace
vectoriel $\Gamma$ de son espace de sections $\H^0(F)$. La dimension
de $\Gamma$ donne l'ordre du syst\`eme coh\'erent. {\`A} l'aide de
r\'esultats de Min He (\cite{He}) sur les espaces de modules de syst\`emes coh\'erents $(\Gamma,F(l))$
d'ordre $1$, dont le \fs sous-jacent est de rang $2$, et de classes de
Chern $c_1=2l, c_2=n+l^2$,
on se ram\`ene au paragraphe 3 de \cite{D1}, pour $n$ compris entre $l(l-1)$
et $(l+1)(l+2)$, \`a l'\'etude de l'espace des
sections d'un  fibr\'e vectoriel $\es^{ld}\R\tens \de^{\tens d}$
sur un ouvert $U$ du sch\'ema de Hilbert $\hil$ des sous-sch\'emas
finis de longueur $m=n+l^2$. Si
$\Xi\subset\hil\times\pp$ est le sous-sch\'ema universel,
 ${\ix}$ est le faisceaux d'id\'eaux associ\'e,
 $pr_1:\hil\times\pp\to \hil$, $pr_2:\hil\times\pp\to \pp$ sont les deux projections, le 
faisceau
alg\'ebrique coh\'erent $\R$ est d\'efini par
$\R=R^1pr_{1*}({\ix}(2l-3))$. Ce faisceau  est localement libre 
en de\-hors du ferm\'e de Brill-Noether $B$ des sch\'emas $Z\in 
\hil$
tels que $h^0(I_Z(2l-3))\ne 0$. On note $U$ l'ouvert  compl\'ementaire de
$B$. La codimension de $B$ est sup\'erieure
ou \'egale \`a $2$, donc les r\'esultats de cohomologie locale  nous
permettent de passer de $\hil$ \`a $U$ pour le calcul d'un espace de sections. Le fibr\'e $\de=\deunu$ est le \fib d\'eterminant sur le sch\'ema de Hilbert
 $\hil$, identifi\'e \`a l'espace de modules  $\em_{(1,0,c_2=n)}$ comme dans le lemme \ref{lemasim}. L'\'enonc\'e pr\'ecis d\'emontr\'e dans \cite{D1} est:
\begin{theor}
\label{th1}
Soit $n$ un entier $\ge 3$. Soit $l$ un entier $>0$ tel que $l(l-1)\le
n<(l+1)(l+2)$. Alors on a un isomorphisme de $\sl(3)$--repr\'esentations
$$\H^0(\em_c,\D^{\tens d})=\H^0(U,\es^{ld}{\mathcal R}\otimes{\mathfrak{d}}^{\tens d}).$$
\end{theor}
On d\'esigne par $E$ l'espace de sections $\H^0(\pp,\O(1))$.
Au paragraphe 4 de \cite{D1} on montre que  $\es^{ld}\R\tens \de^{\tens d}$ admet
sur $U$ une  r\'esolution (*) par un complexe
 $K^{i}=\Lambda^{-i}\es^kE\tens\es^{ld+i}(\vk)\tens\de^{\tens d}$ pour $i=0,\ldots,ld$,
 o\`u $k=2l-3$,  et
 $\vk$ est d\'efini par $\vk=pr_{1*}(\O_{\Xi}\tens pr_2^*(\O(k)))$.

Il est prouv\'e dans \cite{D1} le \th suivant:
\begin{theor}
\label{enunt} 
On a sur $\hil$~:

i) $\H^0(\de^{\tens d})=\es^m(\es^dE)$;

ii) $\H^0(\vk\tens\de^{\tens d})=\es^{k+d}E\tens\es^{m-1}(\es^dE)$;

iii) La ${\sl (3)}$-repr\'esentation  $\H^0(\es^2(\vk)\tens\de^{\tens d})$ est
isomorphe \`a la repr\'esentation
$(\es^{2k+d}E\tens\es^{m-1}(\es^dE))\oplus(\Ker_k\tens\es^{m-2}(\es^dE))$ o\`u
$\Ker_k$ est le noyau de la multiplication
$\es^2(\es^{k+d}E)\to\es^{2k+2d}E$; 

iv)  La ${\sl (3)}$-repr\'esentation $\H^0(\es^3(\vk)\tens\de^{\tens d})$ est
le noyau du morphisme $\alpha$:
$$\alpha:\left[\es^{3k+d}E\tens\es^{m-1}(\es^dE)\right]\oplus\left[\es^{2k+d}E\tens\es^{k+d}E\tens\es^{m-2}(\es^dE)\right]\oplus\left[\es^3(\es^{k+d}E)\tens\es^{m-3}(\es^dE)\right]$$
$$\to\left[\es^{3k+2d-1}E\tens E\tens\es^{m-2}(\es^dE)\right]\oplus\left[\es^{2k+2d}E\tens\es^{k+d}E\tens\es^{m-3}(\es^dE)\right]$$
donn\'e par la matrice:
$$\left(
\begin{array}{ccc}
\tna&\tD&0\\
0&\rho&\tnu\\
\end{array}
\right)
$$
o\`u $\tna$, $\tD$, $\rho$ et $\tnu$ sont des op\'erateurs explicites.

\end{theor}

Le point (i) est une cons\'equence du \th de Kawamata-Viehweg (\cite{C-K-M}),
expliqu\'ee dans le lemme \ref{lemasim}. Le point (ii) correspond au lemme 4.10 de \cite{D1}. Le point (iii) correspond au lemme 4.11 de \cite{D1}. Le point (iv) correspond \`a la proposition 5.13 de \cite{D1}. 

 Le \th \ref{th1} appliqu\'e au
cas $l=1$ nous donne un isomorphisme de $\sl(3)$-repr\'esentations
$$\H^0(\em_c,\D^{\otimes
  d})=\H^0(\Hilb^{n+1}\pp,\es^d\R\otimes\de^{\otimes d}) \mbox{ pour }
  3\le n\le 5.$$
\`A partir de la pr\'esentation (*) de $\R$, appliqu\'ee pour
  $k=2l-3=-1$ et $m=n+l^2=n+1$, on obtient $\vmu\simeq\R$ et
  donc
$$\H^0(\em_c,\D^{\otimes
  d})=\H^0(\Hilb^{n+1}\pp,\es^d(\vmu)\otimes\de^{\otimes d}) \mbox{ pour }
  3\le n\le 5.$$

Pour $d=2$ on applique le \th \ref{enunt} iii) avec $k=2l-3=-1$, $d=2$
et $m=n+l^2=n+1$: $\H^0(\em_c,\D^{\otimes
  2})$ est le noyau du morphisme surjectif
$$\es^n(\es^2E)\oplus\es^2E\otimes\es^{n-1}(\es^2E)\stackrel{(0,id)}{\to}\es^2E\otimes\es^{n-1}(\es^2E).$$
Par suite, la repr\'esentation $\H^0(\em_c,\D^{\otimes
  2})$ est isomorphe \`a $\es^n(\es^2E)$ de dimension $C^n_{n+5}$.
Pour $d=3$ on s'int\'eresse \`a l'espace
$\H^0(\Hilb^{n+1}\pp,\es^3(\vmu)\otimes\de^{\otimes 3})$. Par le \th
\ref{enunt} iv) on est amen\'es \`a \'etudier le noyau du morphisme $\alpha$ (on fait $k=-1$, $d=3$
et $m=n+1$):
\begin{eqnarray*}
  &\alpha:\es^n(\es^3E)\oplus E\otimes
\es^2E\otimes\es^{n-1}(\es^3E)\oplus\es^3(\es^2E)\otimes\es^{n-2}(\es^3E)&\\
&\to E\otimes
\es^2E\otimes\es^{n-1}(\es^3E)\oplus\es^4E\otimes\es^2E\otimes\es^{n-2}(\es^3E).&
\end{eqnarray*}
On montre successivement selon une d\'emarche analogue \`a celle des lemmes 5.16, 5.17, 5.19, 5.20 de 
l'article \cite{D1}, que $\widetilde{D}$ est un isomorphisme, que le
noyau de $\widetilde{\nu}$ est \'egal \`a
$\es^{2,2,2}E\otimes\es^{n-2}(\es^3E)=\comp\otimes\es^{n-2}(\es^3E)$ et son
conoyau \`a $\es^{5,1}E\otimes\es^{n-2}(\es^3E)$, et que
le morphisme de liaison
$(0,\rho):\Ker(\widetilde{\nabla},\widetilde{D})\to\coker\widetilde{\nu}$
est nul. Par cons\'equent, l'\'equivalent de la proposition 5.18 de \cite{D1} nous assure que le noyau de $\alpha$ est isomorphe \`a $\es^n(\es^3E)\oplus
\es^{n-2}(\es^3E)$ de dimension $C^9_{n+9}+C^9_{n+7}$. $\Box$

Ceci conclut la preuve du \th \ref{dual}.

\section{Sections de $\D^{\otimes k}$ pour $n=c_2\le 4$}
\label{3}

Le but de ce paragraphe est de calculer les dimensions des espaces de
sections $\H^0(\em_{(2,0,c_2=n)},\deunu^{\otimes k})$ pour $n\le 4$.

Dans le cas $n=2$, le morphisme de Barth
 nous fournit un \iso entre l'espace de modules $\em_c$ et
$\proj_5$, et le fibr\'e d\'eterminant s'identifie \`a $\O(1)$. Les cas
int\'eressants sont donc $n=3$ et $n=4$. 

On commence par un r\'esultat
 g\'en\'eral sur la fonction $k\mapsto h^0(\emce,\deu^{\otimes k})$ pour
 toute classe $c=(r,c_1,c_2)$ satisfaisant $r>0$ et $\emce$
 non-vide. On rappelle (\S \ref{canodona}) la notation $\gotu=(0,\frac{r}{\delta},-\frac{c_1}{\delta})$ o\`u $\delta={\rm
  pgcd (r,c_1)}$. D'apr\`es le \th \ref{the:pou}, la cohomologie
 sup\'erieure $\H^q(\emce,\deunu^{\otimes k})$ s'annule pour $k\ge
 -3\delta$.

Ce r\'esultat, l'\iso $\omega_{\emce}\simeq\deunu^{\otimes -3\delta}$ et la dualit\'e de Serre nous assurent que
\begin{eqnarray*}
  h^0(\deji)&=&0 {\rm \ \ \ si \ \ \ } j<0\\
h^q(\deji)&=&0 \ \ \ \forall j \ \ \ {\rm \ \ \ si \ } \ \ 0< q< D\\
h^D(\deji)&=&0 {\rm \ \ \ si \ \ \ } j>-3\delta.
\end{eqnarray*}

On note $D=\dim\em_c=1-<c,c^*>$. On note $\D=\D_{\gotu}$. 

\begin{prop}

i) Supposons $d\ge 2$. Pour $k>-3\delta$, la fonction
$$k\mapsto h^0(\em_c,\D^{\otimes k})$$
est un polyn\^ome de degr\'e $D$, de coefficient dominant
$$q_D=\frac{1}{D!}\int_{\em_c}c_1(\D)^D.$$

ii) La s\'erie de Poincar\'e $P(t)$ est de la forme
$$\frac{Q(t)}{(1-t)^{D+1}}$$
o\`u $Q$ est un polyn\^ome de degr\'e  $D+1-3\delta$, \`a
coefficients entiers, tel que $Q(1)=\int_{\em_c}c_1(\D)^D$;

iii) Le polyn\^ome $Q$ satisfait \`a la condition de sym\'etrie
$$t^{D+1-3\delta}Q(\frac{1}{t})=Q(t).$$

\end{prop}

\par{\bf Preuve~:}

La formule de Riemann-Roch pour des vari\'et\'es \'eventuellement
singuli\`eres (\cite{B-F-M}) et le \th \ref{the:pou} donnent
$$\begin{array}{ccccc}
 h^0(\em_c,\D^{\otimes k})&=&\chi(\em_c,\D^{\otimes
  k})&=&\int_{\em_c}ch(\D^{\otimes k})Td(\em_c)\\
&=&\int_{\em_c}e^{kc_1(\D)}Td(\em_c)&=&\sum_{0\le j\le
  D}\frac{k^j}{j!}\int_{\em_c}c_1(\D)^jTd(\em_c).
    \end{array}$$
 Puisque $\D$ est big $q_D>0$, d'o\`u i). 

Posons $a_j=\int_{\em_c}c_1(\D)^jTd(\em_c)$. La s\'erie de Poincar\'e est 
$$\sum_{0\le j\le
  D}a_j(\sum_{k\ge 0}\frac{k^j}{j!}t^k).$$
La somme de la s\'erie $\sum_{k\ge 0}\frac{k^j}{j!}t^k$ (de rayon de
  convergence $1$) est une fonction rationnelle de la forme
  $\frac{Q_j(t)}{(1-t)^{j+1}}$ o\`u les polyn\^omes $Q_j$ sont de degr\'e
  inf\'erieur o\`u \'egal \`a $j$, et $Q_j(1)=1$. Ceci se voit par r\'ecurrence
  sur $j$. Il en r\'esulte que la s\'erie de Poincar\'e est de la forme
  voulue. Puisque $Q(t)=(1-t)^{D+1}P(t)$ et que $P(t)$ est une s\'erie
  formelle \`a coefficients entiers, le calcul de ce produit montre que
  $Q(t)$ est bien \`a coefficients entiers.

La relation classique (voir par exemple \cite{Fulton})
$Td(V^*)ch(\lambda_{-1}(V))=c_{top}(V^*) $ appliqu\'ee au fibr\'e
$V=W\tens\demi$, o\`u $W$ est un espace vectoriel de dimension $m=D+1$,
prouve que 
$$ch(\lambda_{-1}(W\tens \demi))=c_{top}(W^*\tens
\D)Td^{-1}(W^*\tens \D).$$ 
Ici 
$$ \lambda_{-1}(V)=\Lambda^0V-\Lambda^1V+\Lambda^2V-\cdots+(-1)^{\dim
  V}\Lambda^{\dim V}V$$ 
d\'esigne la somme altern\'ee des puissances
ext\'erieures du \fib $V$ et $c_{top}$ d'un \fib vectoriel de rang $r$
d\'esigne la classe de Chern $c_r$ de ce fibr\'e. 

 Mais $V=W\tens\demi$ est un \fib de rang
$D+1$, donc sa classe de Chern maximale $c_{D+1}(W\tens\demi)$
appartient \`a l'espace de cohomologie $\H^{2(D+1)}(\emce)$ qui est nul
en raison de la dimension de $\emce$ ($\dim\emce=D$). Donc
$ch(\lambda_{-1}(W\tens \demi)\tens\deka)=0$ et par la formule de
Riemann-Roch on obtient
$$\chi(\lambda_{-1}(W\tens \demi)\tens\deka)=\int_{\emce}ch(\lambda_{-1}(W\tens \demi)\tens\deka)Td(\emce)=0$$
soit, en tenant compte de
$\Lambda^i(W\tens\demi)=\Lambda^iW\tens\dei$,
$$S_k=\sum_{i=0}^m(-1)^iC^i_m\chi(\dekai)=0$$
quelque soit $k$.

Comme $Q(t)=(1-t)^{D+1}P(t)$, le coefficient $k$-i\`eme de $Q$ s'\'ecrit
\begin{eqnarray*}
  Q_k&=&\left(\begin{array}{c}m\\0\end{array}\right)h^0(\deka)-\left(\begin{array}{c}m\\1\end{array}\right)h^0(\dekaunu)+\cdots+(-1)^m\left(\begin{array}{c}m\\m\end{array}\right)h^0(\dekam)\\
&=&\sum_{i=0}^m(-i)^iC_m^ih^0(\dekai)
\end{eqnarray*}
(en tenant compte du fait que $h^0(\deji)=0$ si $j<0$). Pour
$k<m-3\delta$ on obtient $Q_k=S_k=0$ en raison de l'annulation de la
cohomologie sup\'erieure $\H^q(\deji)$ pour $j<-3\delta, q>0$ et de
$\H^0(\deji)$ pour $j<0$. D'o\`u $\deg\, Q\le m-3\delta$. La condition de
sym\'etrie s'exprime sur la sym\'etrie des coefficients de $Q$~:
$$Q_k=Q_{m-3\delta-k} {\rm \ \ pour \ } 0\le k\le m-3\delta.$$
Mais dans ce cas on a
\begin{eqnarray*}
  Q_k&=&S_k-\sum_{i=k+3\delta}^m(-i)^i\left(\begin{array}{c}m\\i\end{array}\right)\chi(\dekai).
\end{eqnarray*}
La dualit\'e de Serre s'\'ecrit pour $j\le-3\delta$ sous la forme
$\chi(\deji)=(-1)^Dh^0(\dejide)$. En l'appliquant dans la somme
ci-dessus pour $j=k-i\le-3\delta$, et en faisant ensuite la
transformation $j=m-i$ on trouve
\begin{eqnarray*}
  Q_k&=&S_k-(-1)^m\sum_{i=k+3\delta}^m(-i)^{m-i}\left(\begin{array}{c}m\\m-i\end{array}\right)(-1)^Dh^0(\dekaide)\\
&=&S_k+\sum_{j=0}^{m-k-3\delta}(-i)^j\left(\begin{array}{c}m\\j\end{array}\right)h^0(\dekaji)\\
&=&S_k+Q_{m-3\delta-k}=Q_{m-3\delta-k}.
\end{eqnarray*}
En particulier $Q_0=Q_{m-3\delta}=1$ donc le degr\'e de $Q$ est \'egal \`a
$m-3\delta$ exactement. $\Box$

\vs Nous revenons aux cas particuliers qui nous int\'eressent. On prend
$c=(2,0,c_2=n)$ pour $n=3,4$, $\delta=2$ et $\gotu=(0,1,0).$

\par{\bf Preuve du \th \ref{203204}~:}

i) Ici, $D=4c_2-3=9$. La proposition pr\'ec\'edente donne que $P(t)$ s'\'ecrit sous la
forme $\frac{Q(t)}{(1-t)^{10}}$ o\`u $Q$ est de degr\'e $4$, v\'erifie la
condition de sym\'etrie et $Q(1)$ est \'egal \`a $3=9!q_9$ o\`u $q_9$ est un
nombre de Donaldson (\cite{Barth2}). Le calcul de $h^0(\em_c,\D^0)=1$ et
$h^0(\em_c,\D)=10$ (cf. chapitre 2) permet de conclure que $Q(t)=1+t^2+t^4$.

ii) Ce cas est analogue au pr\'ec\'edent seulement il faut faire
intervenir $h^0(\em_c,\D)=15$, et aussi $h^0(\em_c,\D^{\otimes 2})=126$
et $h^0(\em_c,\D^{\otimes 3})=770$ calcul\'es dans \ref{4.0.7}. Ici
$Q(1)=54=13!q_{13}$ (\cite{L-Q}). $\Box$

\begin{rem}
  {\rm Dans le cas $n=3$, puisqu'on a obtenu $h^0(\em_c,\D^{\otimes
  d})$ pour $d=2,3$ (cf. \ref{4.0.7}), on aurait pu d\'eduire la
  valeur du nombre de Donaldson $q_9$.}
\end{rem}

\par{\bf Remerciements~:} Mes remerciements s'adressent \`a J. Le Potier, mon directeur de th\`ese, ainsi qu'\`a N. Dan. Je remercie D. Roessler pour la r\'ef\'erence \cite{Fulton}.

\end{document}